\documentclass[12pt]{amsart} 

\usepackage{amsmath,amssymb,amsthm,amsfonts}
\renewcommand\eqref[1]{(\ref{#1})} 

\usepackage[mathscr]{eucal}
\usepackage{enumitem} 

\setitemize[1]{label=$\bullet$,leftmargin=*}

\newcommand{\abs}[1]{\lvert#1\rvert}
\newcommand{\norm}[1]{\lVert#1\rVert}
\newcommand{\absbig}[1]{\bigl\lvert#1\bigr\rvert}
\newcommand{\normbig}[1]{\bigl\lVert#1\bigr\rVert}
\newcommand{\absBig}[1]{\Bigl\lvert#1\Bigr\rvert}
\newcommand{\normBig}[1]{\Bigl\lVert#1\Bigr\rVert}

\newcommand{\set}[1]{\left\{#1\right\}}
\newcommand{\brac}[1]{(1+\abs{#1})}
\newcommand{\bract}[1]{(1+#1)}

\newcommand\R{{\mathbb R}}
\newcommand\C{{\mathbb C}}
\newcommand\N{{\mathbb N}}

\newcommand\al{\alpha}
\newcommand\be{\beta}
\newcommand\ga{\gamma}
\newcommand\de{\delta}
\newcommand\De{\Delta}
\newcommand\ep{\varepsilon}

\newcommand\la{\lambda}
\newcommand\va{\varphi}
\newcommand\si{\sigma}
\newcommand\Si{\Sigma}
\newcommand\om{\omega}
\newcommand\Om{\Omega}

\newcommand\lap{\Delta}
\newcommand\grad{\nabla}
\newcommand\FT{{\mathscr F}}
\newcommand\pa{\partial}

\newcommand\intbar{-\!\!\!\!\!\!\int}

\renewcommand\S{{\mathcal S}} 

\renewcommand*{\Re}{\operatorname*{Re}}
\renewcommand*{\Im}{\operatorname*{Im}}

\DeclareMathOperator{\meas}{meas}%
\DeclareMathOperator{\supp}{supp}%
\DeclareMathOperator{\codim}{codim}%
\DeclareMathOperator{\Hess}{Hess}%
\DeclareMathOperator{\rank}{rank}%
\DeclareMathOperator{\dist}{dist}%

\newcommand\pat{D_t}%
\newcommand\pax{D_x}%
\newcommand\paxj{D_{x_j}}%
\newcommand\cutoffN{\chi} 
\newcommand\cutoffM{\chi} 
\newcommand\curlyM{\mathcal M}

\numberwithin{equation}{section}
\theoremstyle{plain}
\newtheorem{thm}{Theorem}[section]
\newtheorem{lem}[thm]{Lemma}
\newtheorem{cor}[thm]{Corollary}
\newtheorem{prop}[thm]{Proposition}

\newtheoremstyle{remarks}{}{}{\rmfamily}{}{\bfseries}{}{8pt}
{\thmname{#1}\thmnumber{ #2}:\thmnote{ #3.}}%
\theoremstyle{remarks}
\newtheorem{rem}[thm]{Remark}

\author{Michael Ruzhansky \and  James Smith} 
\title{Global time estimates for solutions to 
equations of dissipative type}
\address{Department of Mathematics \\
Imperial College London \\
180 Queens's Gate \\
London SW7 2AZ \\
United Kingdom }
\email{m.ruzhansky@imperial.ac.uk \\
james.g.smith@imperial.ac.uk}

\begin{document}

\begin{flushleft}
Journees ``Equations aux Derivees Partielles", 
Exp. No. XII, 29 pp., Ecole Polytech., Palaiseau, 2005
\end{flushleft}
\vspace{1.5cm}

\maketitle

\begin{abstract} Global time estimates of $L^p-L^q$ norms of
solutions to general strictly hyperbolic partial 
differential equations are considered. The case of special
interest in this paper
are equations exhibiting the dissipative behaviour.
Results are applied to discuss time decay estimates
for Fokker-Planck equations and for 
wave type equations with negative mass.
\end{abstract}

\section{Introduction}

The paper is devoted to the time
decay of $L^p-L^q$ norms of solutions to 
constant coefficients strictly hyperbolic
equations of general form. 
It is known that such estimates lead
to Strichartz estimates which are a powerful technique 
when dealing
with nonlinear problems. 

We will assume that the principal
part of the equation is strictly hyperbolic. The full equation
may have variable multiplicities because of the lower order
terms. One question of interest is to identify properties
of such equations which determine the time decay rate of
solutions. Another question of interest is what happens when
there are multiple characteristic roots.

Equations of higher orders appear in many applications. 
In particular, they arise as dispersion equations for hyperbolic
systems, for example in the study of
the Fokker-Planck equation and Grad systems in
nonequilibrium thermodynamics. Moreover, in approximations 
of solutions to the Fokker-Planck equation the order of the 
corresponding system tends to infinity. However, it turns
out to still be possible to determine the decay rate of its
solutions. The behaviour exhibited by these examples is 
similar to the behaviour of the dissipative wave equation in
the sense that characteristic roots lie in the 
complex upper half plane and come to the origin 
as single roots and at isolated points. That is why in this
paper we will concentrate on equations of such type in
Theorem \ref{THM:dissipative}, although we will also present
a more general Theorem \ref{THM:overallmainthm}.
Results described here are formulated for scalar equations.
However, they can be easily extended to systems. They also yield
the well-posedness results for semilinear equations.
Details of such analysis will appear elsewhere.

\bigskip
\noindent
{\bf Second order equations.}

\noindent
The study of $L^p-L^q$ decay estimates, or \emph{Strichartz
estimates}, for linear evolution equations began in~1970 when
Robert Strichartz published two papers,~\cite{stri70tran} and
\cite{stri70func}. He proved that if $u=u(x,t)$ satisfies the
Cauchy problem for the
homogeneous linear wave equation
\begin{equation}\label{EQ:waveCauchyprob}
\left\{
\begin{aligned}
&\pa_t^2u(x,t)-\lap_x u(x,t)=0,\quad(x,t)\in\R^n\times(0,\infty)\,,\\
&u(x,0)=\phi(x),\;\pa_t u(x,0)=\psi(x),\quad x\in\R^n\,,
\end{aligned}\right.
\end{equation}
where the initial data $\phi$ and $\psi$ lie in suitable function
spaces such as~$C_0^\infty(\R^n)$, then the \textit{a priori}
estimate
\begin{equation}\label{EQ:strichartzest}
\norm{(u(\cdot,t),u_t(\cdot,t),\nabla_xu(\cdot,t))}_{L^q}\le
C\bract{t}^{-\frac{n-1}{2}\big(\frac{1}{p}-\frac{1}{q}\big)}
\norm{(\nabla_x\phi,\psi)}_{W^{N_p}_p}
\end{equation}
holds when $n\ge2$, $p^{-1}+q^{-1}=1$, $1< p\le2$ and $N_p\ge
n(p^{-1}-q^{-1})$. 
Here $W^{N_p}_p$ stands for the standard Sobolev space
with $N_p$ derivatives over $L^p$.
Using this estimate, Strichartz proved global
existence and uniqueness of solutions to the Cauchy problem for
nonlinear wave equations with suitable (``small'') initial data.
This procedure of proving an \textit{a priori} estimate for a
linear equation and using it, together with local existence of a
nonlinear equation, to prove global existence and uniqueness for a
variety of nonlinear evolution equations is now standard; a
systematic overview, with examples including the equations of
elasticity, Schr\"odinger equations and heat equations, can be
found, for example, in~\cite{rack92}.

There are two main approaches used in order to
prove~\eqref{EQ:strichartzest}; firstly, one may write the
solution to~\eqref{EQ:waveCauchyprob} using the d'Alembert
($n=1$), Poisson ($n=2$) or Kirchhoff ($n=3$) formulae, and their
generalisation to large~$n$,
\begin{equation*}
u(x,t)=\begin{cases}
\begin{aligned} \frac{1}{\prod_{j=1}^{\frac{n-1}{2}}(2j-1)}
\Big[&\pa_t(t^{-1}\pa_t)^{\frac{n-3}{2}}
\Big(t^{n-1}\intbar_{\pa B_t(x)}\phi\,dS\Big)\\
+&(t^{-1}\pa_t)^{\frac{n-3}{2}} \Big(t^{n-1}\intbar_{\pa
B_t(x)}\psi\,dS\Big)\Big]\quad\text{(odd $n\ge3$)}\end{aligned}\\
\begin{aligned} \frac{1}{\prod_{j=1}^{n/2}2j}
&\Big[\pa_t(t^{-1}\pa_t)^{\frac{n-2}{2}} \Big(t^{n}\intbar_{
B_t(x)}\frac{\phi(y)}{\sqrt{t^2-\abs{y-x}^2}}\,dy\Big)\\
+&(t^{-1}\pa_t)^{\frac{n-2}{2}} \Big(t^{n}\intbar_{
B_t(x)}\frac{\psi(y)}{\sqrt{t^2-\abs{y-x}^2}}\,dy\Big)\Big]
\quad\text{(even $n$)}\,,\end{aligned}
\end{cases}
\end{equation*}
(for the derivation of these formulae see, for example,
\cite{evan98}), as is done in~\cite{vonW71} and~\cite{rack92}.
Alternatively, one may write the solution as a sum of Fourier
integral operators:
\begin{equation*}\label{EQ:FMrepforwaveeqn}
u(x,t)=\FT^{-1}
\Big(\frac{e^{it\abs{\xi}}+e^{-it\abs{\xi}}}{2}\,\hat{\phi}(\xi)
+\frac{e^{it\abs{\xi}}-e^{-it\abs{\xi}}}{2\abs{\xi}}\,\hat{\psi}(\xi)
\Big)\,.
\end{equation*}
This is done in~\cite{stri70tran},~\cite{bren75}
and~\cite{pech76}, for example. Using one of these representations
for the solution and techniques from either the theory of Fourier
integral operators (\cite{pech76}), Bessel functions
(\cite{stri70tran}) or standard analysis (\cite{vonW71}), the
estimate~\eqref{EQ:strichartzest} may be obtained. 

Another problem of interest where an $L^p-L^q$ decay estimate for
the linear equation is used to prove existence and uniqueness for
the related nonlinear problem is the Cauchy problem for the
Klein--Gordon equation. Precisely, if $u=u(x,t)$ satisfies the
initial value problem
\begin{equation}\label{EQ:KGCP}
\left\{\begin{aligned} &u_{tt}(x,t)-\lap_xu(x,t)+m^2u(x,t)=0,\quad
(x,t)\in\R^n\times(0,\infty)\,,\\
&u(x,0)=\phi(x),\;u_t(x,0)=\psi(x),\quad x\in\R^n\,,
\end{aligned}\right.
\end{equation}
where $\phi,\psi\in C_0^{\infty}(\R^n)$, say, and $m$ is a
constant (representing a \emph{mass term}), then
\begin{equation}\label{EQ:KGest}
\norm{(u(\cdot,t),u_t(\cdot,t),\nabla_xu(\cdot,t))}_{L^q}\le
C\bract{t}^{-\frac{n}{2}\big(\frac{1}{p}-\frac{1}{q}\big)}
\norm{(\nabla_x\phi,\psi)}_{W^{N_p}_p},
\end{equation}
where $p,q,N_p$ are as before. Comparing~\eqref{EQ:strichartzest}
to~\eqref{EQ:KGest}, we see that the estimate for the solution to
the Klein--Gordon equation decays more rapidly---there is an
improvement in the exponent of the decay function of
$-\frac{1}{2}(\frac{1}{p}-\frac{1}{q})$. The estimate is proved
in~\cite{vonW71},~\cite{pech76} and~\cite{horm97} in different
ways, each suggesting reasons for this improvement:
in~\cite{vonW71}, the function
\begin{equation*}
v=v(x,x_{n+1},t):=e^{-imx_{n+1}}u(x,t)\,,\quad x_{n+1}\in\R\,,
\end{equation*}
is defined; using~\eqref{EQ:KGCP}, it is simple to show that $v$
satisfies the wave equation in $\R^{n+1}$, and thus 
estimate~\eqref{EQ:strichartzest} holds for~$v$, yielding the
desired estimate for~$u$. This is elegant, but cannot easily be
adapted to other situations due to the importance of the
structures of the Klein--Gordon and wave equations for this proof.
In~\cite{pech76} and~\cite{horm97}, a representation of the
solution via Fourier integral operators is used and the stationary
phase method then applied in order to obtain
estimate~\eqref{EQ:KGest}.

A third problem of interest for us is the Cauchy problem for the
dissipative wave equation,
\begin{equation*}
\left\{\begin{aligned}&u_{tt}(x,t)-\lap_x u(x,t)+u_t(x,t)=0\,,
\quad(x,t)\in\R^n\times(0,\infty),
\\ &u(x,0)=\phi(x),\;u_t(x,0)=\psi(x),\quad x\in\R^n\,,
\end{aligned}\right.
\end{equation*}
where $\psi,\phi\in C_0^{\infty}(\R^n)$. In this case,
\begin{equation*}
\norm{\pa_t^r\pa_x^\al u(\cdot,t)}_{L^q}\le
C\bract{t}^{-\frac{n}{2}
(\frac{1}{p}-\frac{1}{q})-r-\frac{\abs{\al}}{2}}
\norm{(\phi,\grad\psi)}_{W_p^{N_p}}\,.
\end{equation*}
This is proved in~\cite{mats76} with a view to showing
well-posedness of related semilinear equations. Once again, this
estimate (for the solution $u(x,t)$ itself) is better than that
for the solution to the wave equation by
$-\frac{1}{2}(\frac{1}{p}-\frac{1}{q})$; there is an even greater
improvement for higher derivatives of the solution. As before, the
proof of this may be done via a representation of the solution
using the Fourier transform:
\begin{equation*}
u(x,t)=\!\!\begin{cases}
\begin{aligned}\FT^{-1}\!\Big(\Big[\frac{e^{-t/2}
\sinh\big(\frac{t}{2}\sqrt{1-4\abs{\xi}^2}\big)}{\sqrt{1-4\abs{\xi}^2}}
+\!e^{-t/2}\!
\cosh\big(\textstyle\frac{t}{2}\sqrt{1-4\abs{\xi}^2}\big)\Big]\!
\hat{\phi}(\xi)\\
+\frac{2e^{-t/2}
\sinh\big(\frac{t}{2}\sqrt{1-4\abs{\xi}^2}\big)}{\sqrt{1-4\abs{\xi}^2}}
\hat{\psi}(\xi)\Big)\,,\quad\abs{\xi}\le1/2,
\end{aligned}\\
\begin{aligned}\FT^{-1}\Big(\Big[\frac{e^{-t/2}
\sin\big(\frac{t}{2}\sqrt{4\abs{\xi}^2-1}\big)}{\sqrt{4\abs{\xi}^2-1}}
+e^{-t/2}
\cos\big(\textstyle\frac{t}{2}\sqrt{4\abs{\xi}^2-1}\big)\Big]
\hat{\phi}(\xi)\\
+\frac{2e^{-t/2}
\sin\big(\textstyle\frac{t}{2}\sqrt{4\abs{\xi}^2-1}\big)}
{\sqrt{4\abs{\xi}^2-1}} \hat{\psi}(\xi)\Big)\,,\quad\abs{\xi}>1/2.
\end{aligned}
\end{cases}
\end{equation*}
Matsumura divides the phase space into the regions where the
solution has different properties and then uses standard
techniques from analysis.

\bigskip
\noindent
{\bf Problem.}

\noindent
It is, therefore, interesting to ask why the addition of lower
order terms improves the rate of decay of the solution to the
equation; furthermore, we would like to understand why the
improvement in the decay is the same for both the addition of a
mass term and for the addition of a dissipative term. In the proof
of each of the estimates (see the papers cited above), the
critical role is played by the characteristic roots of the 
equations. 
In fact, it is the difference in the behaviour of the
characteristic roots of the Klein--Gordon equation and the
dissipative wave equation which yield improvement over the
decay rate for the wave equation.

The aim of this paper is to investigate this phenomenon for
higher order hyperbolic equations and see how lower order terms
affect the rate of decay compared to that for the homogeneous
$m^{\text{th}}$ order equation and the examples above.
Equations of this type appear in many applications. 
In particular, they arise as dispersion equations for hyperbolic
$m\times m$ systems.  The order $m$ may be large,
as in, for example, Grad systems coming from nonequilibrium
thermodynamics, where it corresponds to the number of moments
under consideration. Moreover, in applications to Fokker-Planck
equations describing the distribution of Brownian particles,
the order $m$ corresponds to the Galerkin approximation of 
solutions, so it is increasing to infinity. In all these cases
equations become too large and involved to analyse explicitly,
so we are led to study properties which determine the decay
rate of $L^p-L^q$ estimates in the general form.  

We will consider the Cauchy problem for $m^{{\text{th}}}$ order
constant coefficient linear strictly hyperbolic equation 
of the general form for $u=u(x,t)$:
\begin{equation}\label{EQ:standardCauchyproblem}
\left\{\begin{aligned}& \pat^m
u+\sum_{j=1}^{m}P_{j}(\pax)\pat^{m-j}u+
\sum_{l=0}^{m-1}\sum_{\abs{\al}+r=l}
c_{\al,r}\pax^\al\pat^ru=0,\quad t>0,\\
&\pat^lu(x,0)=f_l(x)\in C_0^{\infty}(\R^n),\quad l=0,\dots,m-1,\;
x\in\R^n\,,
\end{aligned}\right.
\end{equation}
where $P_j(\xi)$ is a constant
coefficient homogeneous polynomial of order~$j$, and the
$c_{\al,r}$ are constants.

We seek \textit{a priori} estimates for the solution to this
problem of the type
\begin{equation}\label{EQ:idealLpLqest}
\norm{\pax^\al\pat^r u(\cdot,t)}_{L^q}\le
K(t)\sum_{l=0}^{m-1}\norm{f_l}_{W^{N_p-l}_p}\,,
\end{equation}
where $1\le p\le2$, $\frac{1}{p}+\frac{1}{q}=1$, $N_p=N_p(\al,r)$
is a constant depending on $p,\al$ and $r$, and $K(t)$ is a
function to be determined.

\bigskip
\noindent
{\bf Homogeneous equations.}

\noindent
The
case where the operator in \eqref{EQ:standardCauchyproblem} is
homogeneous has been studied extensively and provides many 
interesting relations to the geometric properties of characteristics.
In this case we have
\begin{equation}\label{EQ:CPhigherorder}
\left\{\begin{aligned} &L(\pax,\pat)u=0,\quad
(x,t)\in\R^n\times(0,\infty)\\
&\pat^lu(x,0)=f_l(x),\quad l=0,\dots,m-1,\; x\in\R^n\,,
\end{aligned}\right.
\end{equation}
where~$L$ is a homogeneous $m^{\text{th}}$ order constant
coefficient strictly hyperbolic differential operator; the symbol
of~$L$ may be written in the form
\begin{equation*}
L(\tau,\xi)=(\tau-\va_1(\xi))\dots(\tau-\va_m(\xi)),\text{ with }
\va_1(\xi)>\cdots>\va_m(\xi)\quad(\xi\ne0).
\end{equation*}
In a series of papers, \cite{sugi94},~\cite{sugi96}
and~\cite{sugi98}, Sugimoto showed how the geometric
properties of the characteristic roots
$\va_1(\xi),\dots,\va_m(\xi)$ affect the $L^p-L^q$ estimate. To
understand this, let us summarise the method of approach.

Firstly, the solution can be written as the sum of Fourier
multipliers:
\begin{equation*}
u(x,t)=\sum_{l=0}^{m-1}[E_l(t)f_l](x),\quad\text{where }
E_l(t)=\sum_{k=1}^m \FT^{-1}e^{it\va_k(\xi)}a_{k,l}(\xi)\FT
\end{equation*}
and $a_{k,l}(\xi)$ is homogeneous of order $-l$. Now, the problem
of finding an $L^p-L^q$ decay estimate for the solution is reduced
to showing that operators of the form
\begin{equation*}
M_r(D):=\FT^{-1}e^{i\va(\xi)}\abs{\xi}^{-r}\chi(\xi)\FT\,,
\end{equation*}
where $\va(\xi)\in C^\om(\R^n\setminus\set{0})$ is homogeneous of
order~$1$ and $\chi\in C^\infty(\R^n)$ is equal to~$1$ for
large~$\xi$ and zero near the origin, are $L^p-L^q$ bounded for
suitably large~$r\ge l$. In particular, this means that, for
such~$r$,
\begin{equation*}
\norm{E_l(1)f}_{L^q}\le C\norm{f}_{W_p^{r-l}}\,.
\end{equation*}

Indeed, it may be assumed, without loss of generality, that $t=1$
since 
for $t>0$ and $f\in C_0^\infty(\R^n)$,
we easily have 
$$[E_l(t)f](x)=t^l[E_l(1)f(t\cdot)](t^{-1}x)\,.$$
Using this identity gives
\begin{align*}
\norm{E_l(t)f}_{L^q}^q&=
t^{lq}\norm{[E_l(1)f_t](t^{-1}\cdot)}_{L^q}^q
=t^{lq}\int_{\R^n}\abs{[E_l(1)f_t](t^{-1}x)}^q\,dx\\
\stackrel{(x=tx')}{=}&
t^{lq}\int_{\R^n}t^n\abs{[E_l(1)f_t](x')}^q\,dx'
=t^{lq+n}\norm{E_l(1)f_t}_{L^q}^q\,.
\end{align*}
Then, noting that a simple change of variables yields
\begin{equation*}
\norm{f_t}_{W_p^k}^p\le Ct^{kp-n}\norm{f}_{W_p^k}^p\,,
\end{equation*}
we have,
\begin{equation*}
\norm{E_l(t)f}_{L^q} \le Ct^{l+\frac{n}{q}}\norm{f_t}_{W_p^{r-l}}
\le Ct^{r-n(\frac{1}{p}-\frac{1}{q})}\norm{f}_{W_p^{r-l}}\,;
\end{equation*}
hence,
\begin{equation*}
\norm{u(\cdot,t)}_{L^q}\le Ct^{r-n(\frac{1}{p}-\frac{1}{q})}
\sum_{l=0}^{m-1}\norm{f_l}_{W_p^{r-l}}\,.
\end{equation*}

It has long been known that the values of~$r$ for which $M_r(D)$
is $L^p-L^q$ bounded depends on the geometry of the level set
\begin{equation*}
\Si_{\va}=\set{\xi\in\R^n\setminus\set{0}:\va(\xi)=1}\,.
\end{equation*}
In~\cite{litt73},~\cite{bren75} it is shown that if the Gaussian
curvature of $\Si_\va$ is never zero then $M_r(D)$ is $L^p-L^q$
bounded when $r\ge\frac{n+1}{2}\big(\frac{1}{p}-\frac{1}{q}\big)$.
This is extended in~\cite{bren77}, where it is proven that
$M_r(D)$ is $L^p-L^q$ bounded provided
$r\ge\frac{2n-\rho}{2}\big(\frac{1}{p}-\frac{1}{q}\big)$, where
$\rho=\min_{\xi\ne0}\rank\Hess\va(\xi)$.

Sugimoto extended this further in~\cite{sugi94}, where he showed
that if $\Si_\va$ is convex then $M_r(D)$ is $L^p-L^q$ bounded
when $r\ge
\big(n-\frac{n-1}{\ga(\Si)}\big)\big(\frac{1}{p}-\frac{1}{q}\big)$;
here,
\begin{equation*}
\ga(\Si):=\sup_{\si\in\Si}\sup_P\ga(\Si;\si,P)\,,\quad
\Si\subset\R^n\text{ a hypersurface}\,,
\end{equation*}
where $P$ is a plane containing the normal to~$\Si$ at~$\si$ and
$\ga(\Si;\si,P)$ denotes the order of the contact between the line
$T_\si\cap P$, $T_\si$ is the tangent plane at~$\si$, and the
curve $\Si\cap P$. 

In order to apply this result to the solution
of~\eqref{EQ:CPhigherorder}, it is necessary to find a condition
under which the level sets of the characteristic roots are convex.
The following notion is the one that is sufficient.
Let $L=L(\pax,\pat)$ be a homogeneous $m^{\text{th}}$ order
constant coefficient partial differential operator. It is said to
satisfy the \emph{convexity condition} if the Hessian,
$\Hess\varphi_k(\xi)$, corresponding to each of its characteristic
roots $\varphi_1(\xi),\dots,\varphi_m(\xi)$ is semi-definite for
$\xi\ne0$.

It can be shown that if an operator~$L$ does satisfy this
convexity condition, then the above results can be applied to the
solution and thus an estimate of the form \eqref{EQ:idealLpLqest}
holds with
\begin{equation*}
K(t)=\bract{t}^{-\frac{n-1}{\ga}\big(\frac{1}{p}-\frac{1}{q}\big)}\,,
\end{equation*}
where $\gamma=\max_{1\leq k\leq m} \gamma(\Sigma_{\phi_k}).$
We also have $\gamma\leq m$.

Finally, in the case when this convexity condition does not 
hold, it was shown in~\cite{sugi96} and~\cite{sugi98} that, in
general, $M_r(D)$ is $L^p-L^q$ bounded when $r\ge
\big(n-\frac{1}{\ga_0(\Si)}\big)\big(\frac{1}{p}-\frac{1}{q}\big)$,
where
\begin{equation*}
\ga_0(\Si):=\sup_{\si\in\Si}\inf_P\ga(\Si;\si,P)\le \ga(\Si).
\end{equation*}
For $n=2$, $\ga_0(\Si)=\ga(\Si)$, so, the convexity condition may
be lifted in that case. However, in~\cite{sugi96}, examples are
given when $n\ge3$, $p=1,2$ where this lower bound for~$r$ is the
best possible and, thus, the convexity condition is necessary for
the above estimate. It turns out that the case $n\ge3$, $1<p<2$ is
more interesting and is studied in greater depth in~\cite{sugi98},
where microlocal geometric properties must be looked at in order
to obtain an optimal result. It can be noted that in $L^p-L^p$
estimates other geometric properties of phase function and 
wave fronts become important, see the survey \cite{Ruzh} for
more details.

Two remarks are worth making; firstly, the convexity condition
result recovers the Strichartz decay estimate for the wave
equation, since that clearly satisfies such a condition, Secondly,
the convexity condition is an important restriction on the
geometry of the characteristic roots that affects the $L^p-L^q$
decay rate; hence, in the case of an $m^{\text{th}}$ order
operator with lower order terms we must expect some geometrical
conditions on the characteristic roots to obtain decay.

\section{Main Results}

In this paper we will present conditions under which we can
obtain $L^p-L^q$ decay estimates for the general $m^{\text{th}}$
order linear, constant coefficient, strictly hyperbolic Cauchy
problem
\begin{equation}\label{EQ:standardCP(repeat)}
\left\{\begin{aligned}& P(D_t,D_x)\equiv \pat^m
u+\sum_{j=1}^{m}P_{j}(\pax)\pat^{m-j}u+
\sum_{l=0}^{m-1}\sum_{\abs{\al}+r=l}
c_{\al,r}\pax^\al\pat^ru=0,\quad t>0,\\
&\pat^lu(x,0)=f_l(x)\in C_0^{\infty}(\R^n),\quad l=0,\dots,m-1,\;
x\in\R^n\,.
\end{aligned}\right.
\end{equation}
As usual, the strict hyperbolicity means that the principal
symbol of the operator $P(D_t,D_x)$
is strictly hyperbolic, i.e. has
real roots, distinct for $\xi\not=0$. However, characteristic
roots $\tau_1(\xi),\dots,\tau_m(\xi)$ of the full
symbol may have any multiplicities.
Since we are interested in the question of 
how do lower order terms
influence time decay rates, we do not want to worry about the
well-posedness of the Cauchy problem and therefore assume that
the principal part of the operator
is strictly hyperbolic.
Our main Theorem \ref{THM:overallmainthm}
states how different behaviour of the
characteristic roots $\tau_1(\xi),\dots,\tau_m(\xi)$ affect the
rate of decay that can be obtained. 
For now we will assume that symbols of $P_j(D_x)$ are homogeneous
polynomials of order $j$ with constant coefficients. 
However, results extend to the case
when $P_j$ are pseudo-differential operators
of order $j$ and when lower order terms
are pseudo-differential in $D_x$, provided the statement
of Lemma \ref{LEM:multiplerootssetisnice} holds. This
case is essential when considering hyperbolic systems and 
their dispersion equations.

It is natural to impose the condition:
\begin{equation}\label{EQ:imtau>=0}
\Im\tau_k(\xi)\ge0\quad\text{for }k=1,\dots,m\, 
\;\;\textrm{and for all}\;\;\xi\in\R^n;
\end{equation}
this is equivalent to requiring the characteristic polynomial of
the operator to be stable at all points $\xi\in\R^n$, and thus
cannot be lifted, since we can not expect any time decay if
this condition fails. Coefficients of equation \eqref{EQ:standardCP(repeat)}
are allowed to be complex as long as condition \eqref{EQ:imtau>=0}
holds. Of course, because of the strict hyperbolicity,
coefficients of the principal part are real. 

Also, it is sensible to divide the
considerations of how characteristic roots behave into two parts:
their behaviour for large values of~$\abs{\xi}$  and for bounded
values of~$\abs{\xi}$. These two cases are then subdivided
further; in particular the following are the key properties to
consider:
\begin{itemize}
\item multiplicities of roots (this only occurs in the case of
bounded~$\abs{\xi}$);
\item whether roots lie on the real axis or are separated
from it;
\item behaviour as $\abs{\xi}\to\infty$ (only in the case of
large~$\abs{\xi}$);
\item how roots meet the real axis (if they do);
\item properties of the Hessian of the root, $\Hess\tau_k(\xi)$;
\item a convexity-type condition, as in the case of homogeneous
roots.
\end{itemize}

Some definitions will be needed for the main theorem.
Given a smooth function $\tau:\R^n\to\R$ and $\la\in\R$, set
\begin{equation*}
\Si_\la\equiv\Si_\la(\tau):=\set{\xi\in\R^n:\tau(\xi)=\la}\,.
\end{equation*}
In the case where $\tau(\xi)$ is homogeneous of order one, write
$\Si_\tau:=\Si_1(\tau)$---for such $\tau$, we then have
$\Si_\la(\tau)=\la\Si_\tau$.
Also, a smooth function $\tau:\R^n\to\R$ will be 
said to satisfy the
\emph{convexity condition} if~$\Si_\la$ is convex for each
$\la\in\R$. Note that the empty set is considered to be convex.
Finally, we will use notation in the introduction 
for the \emph{maximal orders of
contact} of a hypersurface,
$\ga(\Si)$ and $\ga_0(\Sigma)$.
We note that if $p(\xi)$ is a polynomial of order $m$ and
$\Si=\set{\xi\in\R^n:p(\xi)=0}$ is compact then
$\ga_0(\Si)\le\ga(\Si)\le m$; this is useful when applying the
result below to hyperbolic differential equations and is proved in
\cite{sugi96}.

Now we may state the main theorem:
\begin{thm}\label{THM:overallmainthm}
Suppose $u=u(x,t)$ satisfies the $m^{\text{th}}$ order
linear\textup{,} constant coefficient\textup{,} strictly
hyperbolic Cauchy problem~\eqref{EQ:standardCP(repeat)}. Denote
the characteristic roots of the operator by
$\tau_1(\xi),\dots,\tau_m(\xi)$ and assume
that~\eqref{EQ:imtau>=0} holds.

We introduce two functions\textup{,} $K^{(\text{l})}(t)$ and
$K^{(\text{b})}(t)$, which take values as follows\textup{:}
\begin{enumerate}[leftmargin=*,label=\textup{\Roman*.},ref=\Roman*]
\item\label{RESULT:mainthmlargexi}
Consider the behaviour of~each characteristic root\textup{,}
$\tau_k(\xi)$\textup{,} in the region $\abs{\xi}\ge N$\textup{,}
where~$N$ is some large number. The following
table gives values for the function $K_k^{(\text{l})}(t)$
corresponding to possible properties of~$\tau_k(\xi)$\textup{;} if
$\tau_k(\xi)$ satisfies more than one\textup{,} then take
$K_k^{(\text{l})}(t)$ to be function that decays the slowest as
$t\to\infty$.
\end{enumerate}
\begin{center}\begin{upshape}
\begin{tabular}[4]{|c|c|c|}\hline
Location of $\tau_k(\xi)$ & Additional Property &
$K_k^{\text(l)}(t)$\\
\hline\hline%
away from real axis && $e^{-\de t}$, some $\de>0$\\\hline
&$\det\Hess\tau_k(\xi)\ne0$ &
$\bract{t}^{-\frac{n}{2}(\frac{1}{p}-\frac{1}{q})}$\\
on real axis &$\rank\Hess\tau_k(\xi)=n-1$ &
$\bract{t}^{-\frac{n-1}{2}(\frac{1}{p}-\frac{1}{q})}$\\
& convexity condition, $\ga$ &
$\bract{t}^{-\frac{n-1}{\ga}(\frac{1}{p}-\frac{1}{q})}$\\ & no
convexity condition, $\ga_0$ &
$\bract{t}^{-\frac{1}{\ga_0}(\frac{1}{p}-\frac{1}{q})}$\\
\hline &$\det\Hess\tau_k(\xi)\ne0$ &
$\bract{t}^{-\frac{n}{2}(\frac{1}{p}-\frac{1}{q})}$\\
asymptotic to real axis &$\rank\Hess\tau_k(\xi)=n-1$ &
$\bract{t}^{-\frac{n-1}{2}(\frac{1}{p}-\frac{1}{q})}$\\
& no convexity condition, $\ga_0$ &
$\bract{t}^{-\frac{1}{\ga_0}(\frac{1}{p}-\frac{1}{q})}$\\ \hline
\end{tabular}\end{upshape}
\end{center}
Then take
$K^{(\text{l})}(t)=\max_{k=1\,\dots,m}K_k^{\text(l)}(t)$.

\begin{enumerate}[leftmargin=*,label=\textup{\Roman*.},resume,ref=\Roman*]
\item\label{RESULT:mainthmbddxi}
Consider the behaviour of the characteristic roots in the bounded
region~$\abs{\xi}\le N$\textup{;} again\textup{,} take
$K^{(\text{b})}(t)$ to be the maximum \textup{(}slowest
decaying\textup{)} function for which there are roots satisfying
the conditions in the following table\textup{:}
\end{enumerate}
\begin{center}\begin{upshape}
\begin{tabular}[4]{|c|c|c|}\hline
Location of Root(s)& Properties & $K^{(\text{b})}(t)$\\
\hline\hline%
away from axis & no multiplicities &$e^{-\de t}$, some $\de>0$\\
&$L$ roots coinciding & $\bract{t}^{L-1} e^{-\de t}$\\
\hline on axis,&$\det\Hess\tau_k(\xi)\ne0$ &
$\bract{t}^{-\frac{n}{2}(\frac{1}{p}-\frac{1}{q})}$\\
no multiplicities & convexity condition, $\ga$ &
$\bract{t}^{-\frac{n-1}{\ga}(\frac{1}{p}-\frac{1}{q})}$\\
& no convexity condition, $\ga_0$ &
$\bract{t}^{-\frac{1}{\ga_0}(\frac{1}{p}-\frac{1}{q})}$\\
\hline meeting axis & $L$ roots coincide&\\
with finite order $s$& on set of codimension $\ell$&
$\bract{t}^{L-1-\frac{\ell}{s}(\frac{1}{p}-\frac{1}{q})}$\\
\hline
\end{tabular}\end{upshape}
\end{center}

Then, with
$K(t)=\max\big(K^{\text{(b)}}(t),K^{\text{(l)}}(t)\big)$\textup{,}
the following estimate holds\textup{:}
\begin{equation*}
\norm{\pax^\al\pat^r u(\cdot,t)}_{L^q}\le K(t)
\sum_{l=0}^{m-1}\norm{f_l}_{W^{N_p-l}_p}\,,
\end{equation*}
where $1< p\le2$, $\frac{1}{p}+\frac{1}{q}=1$, and
$N_p=N_p(\al,r)$ is a constant depending on~$p,\al$ and~$r$.
\end{thm}
Let us make a number of remarks on how to understand this theorem.
Since the decay rate does depend on the behaviour of characteristic
roots at different points, we single out properties which determine
this decay rate. Since the same characteristic root, say $\tau_k$,
may exhibit different properties at different points, we look
at the corresponding rates 
$K^{\text{(b)}}(t),K^{\text{(l)}}(t)$ under each possible condition
and then take the slowest one for the final answer. It also means
that if we microlocalise in a region where only one of these
properties holds, we can get 
the decay rate straight from the table for
the corresponding solution. 
In some cases, especially when roots do not lie on the axis
for large $\xi$, the result may be extended to $p=1$.

In Part I of the statement, it can be
shown by the perturbation arguments that only three
cases are possible for large $\xi$, namely, the characteristic root
may be uniformly separated from the real axis, it may lie on the
axis, or it may converge to the real axis at infinity. If,
for example, the
root lies on the axis and, in addition, it satisfies the convexity
condition with index $\gamma$, we get the corresponding decay rate
$K^{\text{(l)}}(t)=\bract{t}^{-\frac{n-1}{\ga}(\frac{1}{p}-\frac{1}{q})}$.
Indices $\gamma$ and $\gamma_0$ in the tables are defined
as the maximum of the corresponding indices $\gamma(\Sigma_\lambda)$
and $\gamma(\Sigma_\lambda)$, where 
$\Sigma_\lambda=\{\xi:\tau_k(\xi)=\lambda\}$, over all $k$ and over
all $\lambda$, for which $\xi$ lies in the corresponding zone.

The statement in Part II is more involved since we may have multiple
roots intersecting on rather irregular sets. 
The number $L$ of coinciding roots corresponds to the number of
roots which actually contribute to the loss of regularity.
For example, operator $(\partial_t^2-\Delta)(\partial_t^2-2\Delta)$
would have $L=2$ for both pairs of roots intersecting at the
origin.  Meeting the axis with finite order $s$ 
means that we have the estimate 
\begin{equation}\label{EQ:disttau}
\dist(\xi,Z_k)^s\leq c|\Im\tau_k(\xi)|
\end{equation}
for all the intersecting roots, where $Z_k=\{\xi: \Im\tau_k(\xi)=0\}.$
In Part II of Theorem \ref{THM:overallmainthm}, 
the condition that 
$L$ roots meet the axis with finite order $s$ on a set of codimension
$\ell$ means that all these estimates hold
and that there is a (regular) set $Z$ 
of codimension $\ell$ such that
$Z_k\subset Z$ for all corresponding $k$. 
In Theorem \ref{THM:dissipative} we will discuss the special case
of a single root $\tau_k$
meeting the axis at a point $\xi_0$ with order $s$,
which means that $\Im\tau_k(\xi_0)=0$ and 
that we have the estimate 
$\abs{\xi-\xi_0}^s\le c\abs{\Im\tau_k(\xi)}$.
In fact, under certain conditions an improvement in this part
of the estimates is possible, see Theorem \ref{THM:dissipative}
and Remark \ref{REM:dismore}. 

In addition to the theorem, if we have $L$ multiple roots 
which coincide on the real axis on a set $S$ of codimension $\ell$,
we have an estimate  
\begin{equation}\label{EQ:around}
|u(t,x)|\leq C(1+t)^{L-1-\ell}
\sum_{l=0}^{m-1}\norm{f_l}_{L^1},
\end{equation}
if we cut off the Fourier transform of the Cauchy data
to the $\epsilon$-neighbourhood $S^\epsilon$ of $S$
with $\epsilon=1/t$. 
Here we may  relax the definition of the intersection
above and say that if $L$ roots coincide
on a set $S$, then 
they coincide on a set of codimension $\ell$
if the measure of the $\epsilon$-neighborhood $S^\epsilon$ of $S$
satisfies $|S^\epsilon|\leq C\epsilon^{\ell}$ for small $\epsilon>0$;
here $S^\epsilon=\{\xi\in\R^n: \dist (\xi,S)\leq\epsilon\}.$
The estimate \eqref{EQ:around}
follows from the procedure
described below of the resolution
of multiple roots. 
We can then combine this with the remaining cases outside
of this neighborhood, where it is possible to establish
decay by different arguments. In particular, this is the case
of homogeneous equations with roots intersecting at the
origin. However, one sometimes needs to introduce special
norms to handle $L^2$-estimates around the multiplicities.
Details of this will appear elsewhere.
Finally, in the case of a simple root we may set $L=1$,
and $\ell=n$, if it meets the axis at a point.

Theorem \ref{THM:overallmainthm} allows a microlocalisation and
estimates for the corresponding oscillatory integrals.
In fact, Theorem \ref{THM:overallmainthm} follows from its 
microlocal version in regions where characteristic roots are
simple. In regions with multiple roots one requires additional
arguments resolving the singularities caused by multiple roots
followed by estimates for relevant pieces of the solution.
A microlocal version of the theorem leads to better estimates
since in Theorem \ref{THM:overallmainthm} we take the slowest
among all microlocal decay rates. Such a microlocal version
and the full proof of Theorem
\ref{THM:overallmainthm} will appear elsewhere.

For our applications in Section~\ref{SEC:applications}, 
we only need this result in the special case where 
characteristic roots meet the real axis
with finite order; 
therefore, we shall state and outline the prove the theorem 
in this special case. 

\begin{thm}
\label{THM:dissipative}
Consider the $m^{\text{th}}$ order strictly hyperbolic Cauchy
problem~\eqref{EQ:standardCP(repeat)} for operator
$P(D_t,D_x)$, with initial data $f_j\in
L^p\cap W^{\left[\frac{n}{2}\right]+1+\abs{\al}-j+r}_2$ for
$j=0,\dots,m-1$, where $1\le p\le2$ and $2\leq q\leq\infty$
are such that $\frac{1}{p}+\frac{1}{q}=1$.
Assume that the
characteristic roots $\tau_1(\xi),\dots,\tau_m(\xi)$ 
of $P(\tau,\xi)=0$
satisfy~\eqref{EQ:imtau>=0}
and the following conditions:
\begin{enumerate}[leftmargin=*,label=\textup{(H\arabic*)}]
\item there is some $\epsilon>0$ such that 
for all $k\in\set{1,\dots,m}$ we have
$$\liminf_{\abs{\xi}\to\infty}\Im\tau_k(\xi)\ge\ep\,;$$
\item for each
$\xi_0\in\R^n$ there is at most one $k$
for which $\Im\tau_k(\xi_0)=0$ and there exists
a constant $c>0$ such that
\begin{equation*}
\abs{\xi-\xi_0}^s\le c\abs{\Im\tau_k(\xi)},
\end{equation*}
for $\xi$ in some neighbourhood of $\xi_0$.
\end{enumerate}
Then the solution $u=u(x,t)$ to 
Cauchy problem \eqref{EQ:standardCP(repeat)} satisfies 
the estimate
%
\begin{equation*}\label{EQ:LpLqest}
\norm{\pat^r\pax^\al u(\cdot,t)}_{L^q}\le
C_{\al,r}\bract{t}^{-\frac{n}{s}
(\frac1p-\frac1q)}\sum_{j=0}^{m-1}
\norm{f_j}_{W^{\frac{2-p}{p}\big(\left[\frac{n}{2}\right]+1\big)+
\abs{\al}+r-j}_p}.
\end{equation*}
\end{thm}
Essentially, this theorem is a special case of Theorem
\ref{THM:overallmainthm}, where we get the exponential decay
from Part I, exponential decay from multiple roots away from
the real axis in Part II, as well as the last line of the
table in Part II with $L=1$ and $l=n$, since we have only
a single root coming to the axis. The main problem here is 
the possible appearance of multiple roots in the complex
upper half plane. If several roots meet on the axis,
the decay is then given in Part II of Theorem
\ref{THM:overallmainthm}, where we observe the appearance
of the extra power $t^{L-1}$ compared to Theorem
\ref{THM:dissipative}. If roots come to the axis on a set
of other codimension $\ell$, the order should change
according to Theorem \ref{THM:overallmainthm}.
If conditions of Theorem \ref{THM:dissipative} hold only with
$\xi_0$ we will call the polynomial $P(\tau,\xi)$
strongly stable. Such polynomials will be discussed in more
detail in Section \ref{SEC:applications}.

\begin{rem}\label{REM:dismore}
The order of time decay in Theorem \ref{THM:dissipative}
may be improved in the following cases.
If $\Im\tau_k(\xi_0)=0$ in (H2) 
implies that $\xi_0=0$, then we actually get
\begin{equation*}
\normBig{D^r_tD^\al_x
u(\cdot,t)}_{L^q(\R^n_x)}\le
C\bract{t}^{-\frac{n+\abs{\al}}{s}\big(\frac{1}{p}-\frac{1}{q}\big)}
\sum_{j=0}^{m-1}
\norm{f_j}_{W^{\frac{2-p}{p}\big(\left[\frac{n}{2}\right]+1\big)+
\abs{\al}+r-j}_p}
\,.
\end{equation*}
Now, assume that for all $\xi_0$ in (H2) we also have
the estimate
\begin{equation}\label{EQ:imest}
c_0\abs{\xi-\xi_0}^{s}\le c\abs{\Im\tau_k(\xi)}\leq
c_1\abs{\xi-\xi_0}^{s_1}.
\end{equation}
with some constants $c_0,c_1>0$. 

If $\Im\tau_k(\xi_0)=0$ in (H2) 
implies that $\Re\tau_k(\xi_0)=0$, 
then we actually get
\begin{equation*}
\normBig{D^r_tD^\al_x
u(\cdot,t)}_{L^q(\R^n_x)} \le
C\bract{t}^{-\big(\frac{n+r s_1}{s}\big)\big(\frac{1}{p}-\frac{1}{q}\big)}
\sum_{j=0}^{m-1}
\norm{f_j}_{W^{\frac{2-p}{p}\big(\left[\frac{n}{2}\right]+1\big)+
\abs{\al}+r-j}_p}
\,.
\end{equation*}

And finally, assume that for all $\xi_0$ such that
$\Im\tau_k(\xi_0)=0$ in (H2), we also have
$\xi_0=0$ and $\Re\tau_k(\xi_0)=0$.
Then we actually get
\begin{equation*}
\normBig{D^r_tD^\al_x
u(\cdot,t)}_{L^q(\R^n_x)}\le
C\bract{t}^{-\frac{n+\abs{\al}+r s_1}{s}\big(\frac{1}{p}-\frac{1}{q}\big)}
\sum_{j=0}^{m-1}
\norm{f_j}_{W^{\frac{2-p}{p}\big(\left[\frac{n}{2}\right]+1\big)+
\abs{\al}+r-j}_p}
\,.
\end{equation*}
\end{rem}
The proof of this is based on Remark \ref{rem:mats}.
In particular, these estimates cover the case of dissipative
wave equation and applications in Section 
\ref{SEC:applications}. A similar remark can be made for Theorem
\ref{THM:overallmainthm}, where we also get the corresponding
improvements.

\section{Outline of the proof}

Here we will outline the proof of Theorem \ref{THM:dissipative}.
For large frequencies we have simple roots separated from the
real axis so we can expect exponential decay in time
there. For small frequencies, while separated from the real
axis, we may have multiple roots, which may intersect on a 
rather irregular set. We will cut off around this set and show
that we can get additional polynomial growth in time dependent
on the ``dimension'' of this set, which is matched against 
exponential decay. Technically we have to establish a number
of additional estimates on the solution in this case since the
usual solution representation blows up around points of
multiplicity.
Finally, we can show the polynomial decay in time when 
characteristic roots approach the real axis. 

\subsection{Some properties of hyperbolic polynomials}

Here we will describe some useful properties of hyperbolic
polynomials.
Let $L=L(\pax,\pat)$ be a linear $m^{\text{th}}$ order constant
coefficient partial differential operator. Then each of the
characteristic roots of $L$\textup{,} denoted
$\tau_1(\xi),\dots,\tau_m(\xi)$\textup{,} is continuous in
$\R^n$\textup{;} furthermore\textup{,} for each
$k=1,\dots,m$\textup{,} the characteristic root $\tau_k(\xi)$ is
analytic in
\begin{equation*}\set{\xi\in\R^n:\tau_k(\xi)\ne
\tau_l(\xi)\,\forall\,l\ne k}\,.\end{equation*}
Let now $L=L(\pax,\pat)$ be a linear $m^{\text{th}}$ order constant
coefficient strictly hyperbolic partial differential operator. Then
there exists a constant $N$ such that\textup{,} 
the characteristic roots $\tau_1(\xi),\dots,\tau_m(\xi)$ of $L$
are pairwise distinct for $|\xi|\geq N$.
We also have the following symbolic properties of characteristic
roots:

\begin{prop}\label{PROP:perturbationresults} 
Let $L=L(\pax,\pat)$ be a
linear $m^{\text{th}}$ order constant coefficient 
hyperbolic partial
differential operator with characteristic roots
$\tau_1(\xi)$, $\dots$, $\tau_m(\xi)$\textup{;} then
\begin{enumerate}[label=\upshape{\Roman*.},ref={Part~\Roman*}]
\item\label{LEM:orderroot} for each $k=1,\dots,m$\textup{,} there
exists a constant $C>0$ such that
\begin{equation*}\label{EQ:tauisoderxi}
\abs{\tau_k(\xi)}\le C\brac{\xi}\quad\text{for all }\xi\in\R^n\,.
\end{equation*}
\end{enumerate}
Suppose that the maximum order of
the lower order terms is $0\le K\le m-1$.
Furthermore\textup{,} assume that~$L$ is \emph{strictly}
hyperbolic\textup{,} and denote the roots of the principal part
$L_m(\xi,\tau)$ by $\varphi_1(\xi),\dots,\varphi_m(\xi)$.
Then we have the following\textup{:}

\begin{enumerate}[label=\upshape{\Roman*.},ref={Part~\Roman*},resume]
\item\label{LEM:tau-vabounds}  For each
$\tau_k(\xi)$\textup{,} $k=1,\dots,m$\textup{,} there exists a
corresponding root of the principal symbol $\va_k(\xi)$
\textup{(}possibly after reordering\textup{)} such that
\begin{equation*}\label{EQ:tau-vaboundkthorder}
\abs{\tau_k(\xi)-\va_k(\xi)}\le C\brac{\xi}^{K+1-m} \quad\text{for
all }\xi\in\R^n\,.
\end{equation*}

\item\label{PROP:boundsonderivsoftau} There exists $N>0$
such that\textup{}, for each characteristic root of~$L$ and for each
multi-index $\al$\textup{,} we can find constants $C=C_{k,\al}>0$
such that
\begin{equation*}\label{EQ:rootisassymbol}
\absbig{\pa^\al_\xi\tau_k(\xi)}\le
C\abs{\xi}^{1-\abs{\al}}\,\quad\text{for all }\abs{\xi}\ge N\,,
\end{equation*}

\item \label{PROP:derivoftau-derivofphi} There exists $N>0$ such
that\textup{,} for each $\tau_k(\xi)$ a corresponding root of the
principal symbol~$\va_k(\xi)$ can be found \textup{(}possibly after
reordering\textup{)} which satisfies\textup{,} for each multi-index
$\al$ and $k=1,\dots,m$\textup{,}
\begin{equation*}\label{EQ:derivsoftau-phiforfewerlot}
\absbig{\pa^\al_\xi\tau_k(\xi)-\pa_\xi^\al\va_k(\xi)} \le
C\abs{\xi}^{K+1-m-\abs{\al}}\quad\text{for all }\abs{\xi}\ge N
\end{equation*}
for each multi-index $\al$ and $k=1,\dots,m$.
\end{enumerate}
\end{prop}

\subsection{Representation of the solution}

Recall that we begin with the Cauchy problem with solution
$u=u(x,t)$:
\begin{equation}\label{EQ:standardCP(ch4)}
\left\{\begin{aligned}& \pat^m
u+\sum_{j=1}^{m}P_{j}(\pax)\pat^{m-j}u+
\sum_{l=0}^{m-1}\sum_{\abs{\al}+r=l}
c_{\al,r}\pax^\al\pat^ru=0,\quad t>0,\\
&\pat^lu(x,0)=f_l(x)\in C_0^{\infty}(\R^n),\quad l=0,\dots,m-1,\;
x\in\R^n\,,
\end{aligned}\right.
\end{equation}
where symbol $P_j(\xi)$ of $P_j(D)$  is a constant
coefficient homogeneous polynomial of order~$j$, and
the~$c_{\al,r}$ are constants.

Applying the partial Fourier transform with respect to~$x$ yields an
ordinary differential equation for $\hat{u}=
\hat{u}(\xi,t):=\int_{\R^n}e^{-ix\cdot\xi}u(x,t)\,dx$:
\begin{subequations}
\begin{align}\label{EQ:ftcauchyprob}
\pat^m\hat{u}+\sum_{j=1}^mP_j(\xi)\pat^{m-j}\hat{u}+
\sum_{l=0}^{m-1}\sum_{\abs{\al}+r=l}c_{\al,r}
\xi^\al\pat^r\hat{u}=0\,,
\\
\pat^l\hat{u}(\xi,0)=\hat{f_l}(\xi),\quad
l=0,\dots,m-1,\label{EQ:ftcauchydata}
\end{align}
where $(\xi,t)\in\R^n\times[0,\infty)$. Let $E_j=E_j(\xi,t)$,
$j=0,\dots,m-1$, be the solutions to~\eqref{EQ:ftcauchyprob} with
initial data
\begin{equation}\label{EQ:initialdataforEj}
\pat^lE_j(\xi,0)=\begin{cases}1\quad\text{ if }l=j,\\
0\quad\text{ if }l\ne j.\end{cases}
\end{equation}
\end{subequations}
Then the solution~$u$ of~\eqref{EQ:standardCP(ch4)} can be written
in the form
\begin{equation}\label{EQ:repforu}
u(x,t)=\sum_{j=0}^{m-1}(\FT^{-1}E_j\FT f_j)(x,t),
\end{equation}
where $\FT$ and $\FT^{-1}$ represent the partial Fourier transform
with respect to $x$ and its inverse respectively.

Now, as~\eqref{EQ:ftcauchyprob},~\eqref{EQ:initialdataforEj} is
the Cauchy problem for a linear ordinary differential equation, we
can write, denoting the characteristic roots
of~\eqref{EQ:standardCP(ch4)} by $\tau_1(\xi),\dots,\tau_m(\xi)$,
\begin{equation*}
E_j(\xi,t)=\sum_{k=1}^mA^k_j(\xi,t)e^{i\tau_k(\xi)t}
\end{equation*}
where $A^k_j(\xi,t)$ are polynomials in~$t$ whose coefficients
depend on~$\xi$. Moreover, for each $k=1,\dots,m$ and
$j=0\dots,m-1$, the $A_j^k(\xi,t)$ are independent of~$t$ at
points of the (open) set $\set{\xi\in\R^n:\tau_k(\xi)\ne
\tau_l(\xi)\,\forall\,l\ne k}$; when this is the case, we write
$A_j^k(\xi,t)\equiv A_j^k(\xi)$.  For $A_j^k(\xi)$,
we have the following properties:
\begin{lem}\label{LEM:ordercoeff}
Suppose $\xi\in S_k:=\set{\xi\in\R^n:\tau_k(\xi)\ne
\tau_l(\xi)\,\forall\,l\ne k}$\textup{;} then we have the
following formula\textup{:}
\begin{equation}\label{EQ:Ajkformula}
A_j^k(\xi)=\frac{(-1)^j\displaystyle\sideset{}{^k}\sum_{1\le
s_1<\dots<s_{m-j-1}\le
m}\prod_{q=1}^{m-j-1}\tau_{s_q}(\xi)}{\displaystyle\prod_{l=1,l\ne
k}^m(\tau_l(\xi)-\tau_k(\xi))}\;,
\end{equation}
where $\sum^k$ means sum over the range indicated excluding $k$.
Furthermore\textup{,} we have\textup{,} for each $j=0,\dots,m-1$
and $k=1,\dots,m$\textup{,}
\begin{enumerate}[label=\textup{(}\roman*\textup{)},leftmargin=*]
\item $A_j^k(\xi)$ is smooth in~$S_k$\textup{;}
\item $A_j^k(\xi)=O(\abs{\xi}^{-j})$ as $\abs{\xi}\to\infty$.
\end{enumerate}
\end{lem}
\begin{proof}
The representation~\eqref{EQ:Ajkformula} follows from Cramer's rule
(and is done explicitly in~\cite{klin67}): $A_j^k(\xi)=\frac{\det
V_j^k}{\det V}$, where $V:=\big(\tau_i^{l-1}(\xi)\big)_{i,l=1}^m$
is the Vandermonde matrix and~$V_j^k$ is the matrix obtained by
taking~$V$ and replacing the $k^{\text{th}}$ column by
$(\underbrace{0\ \dots\ 0\ 1}_j\ 0\ \dots\ 0)^{\mathrm T}$.

Smoothness of $A_j^k(\xi)$ in $S_k$ is obvious 
and the asymptotic behaviour is
a consequence of \ref{LEM:orderroot} of
Proposition~\ref{PROP:perturbationresults}
since~\eqref{EQ:Ajkformula} holds for all $\abs{\xi}>N$.
\end{proof}

In view of Lemma \ref{LEM:ordercoeff},
choose $N_1>0$ so that the $\tau_k(\xi)$, $k=1,\dots,n$, are
distinct for $\abs{\xi}>N_1$. Also, choose~$N_2>0$ so that all
points at which any of the roots, $\tau_k(\xi)$, meet the real
axis---i.e.\ points $\xi\in\R^n$ such that, for all $\ep>0$, there
exist $\xi_1,\xi_2\in B_\ep(\xi)$ with $\Im\tau_k(\xi_1)=0$ and
$\Im\tau_k(\xi_2)\ne0$---lie in~$B_{N_2}(0)$. Set
$N=\max(N_1,N_2)$.

Let $\cutoffN(\xi)=\cutoffN_N(\xi)\in C_0^\infty(\R^n)$,
$0\le\cutoffN(\xi)\le1$, be a cut-off function that is
identically~$1$ for $\abs{\xi}<N$ and identically zero for
$\abs{\xi}>2N$. Then \eqref{EQ:repforu} can be rewritten as:
\begin{equation}\label{EQ:divideduprepforu}
u(x,t)=\sum_{j=0}^{m-1}\FT^{-1}(E_j\cutoffN\FT f_j)(x,t)
+\sum_{j=0}^{m-1}\FT^{-1}(E_j(1-\cutoffN)\FT f_j)(x,t)\,.
\end{equation}

\subsection{Large $\abs\xi$} 
The second term of~\eqref{EQ:divideduprepforu} is the most straightforward
to study: by the choice of~$N$,
\begin{equation*}
E_j(\xi,t)(1-\cutoffN)(\xi)=\sum_{k=1}^mA_j^k(\xi)(1-\cutoffN)(\xi)
e^{i\tau_k(\xi)t}\,;
\end{equation*}
therefore, since each summand is smooth in~$\R^n$,
\begin{multline*}
\sum_{j=0}^{m-1}\FT^{-1}(E_j(1-\cutoffN)\FT
f_j)(x,t)\\=\frac{1}{(2\pi)^{n}}\sum_{j=0}^{m-1}\sum_{k=1}^m
\int_{\R^n} e^{i(x\cdot\xi+\tau_k(\xi)t)}
A_j^k(\xi)(1-\cutoffN)(\xi)\hat{f}_j(\xi)\,d\xi\,.
\end{multline*}
Note that, unlike in the case of homogeneous strictly hyperbolic
equations we may not assume that $t=1$. Each of these integrals may be studied
separately. Indeed, we have the following result:
\begin{prop}\label{PROP:rootsawayfromaxis}
Let $\tau:U\to\C$ be a smooth function, $U\subset\R^n$ open, and
$a_j\in\S^{-j}_{1,0}(U)$. Assume\textup{:}
\begin{enumerate}[leftmargin=*,label=\textup{(\roman*)}]
\item there exists $\de>0$ such that $\Im\tau(\xi)\ge\de$ for all
$\xi\in U$\textup{;}
\item $\abs{\tau(\xi)}\le C\brac{\xi}$ for all $\xi\in U$.
\end{enumerate}
Then\textup{,}
\begin{gather*}
\normBig{\int_{U}e^{i(x\cdot\xi+\tau(\xi)t)}a_j(\xi)
\xi^{\al}\tau(\xi)^r \hat{f}(\xi)\,d\xi}_{L^\infty(\R^n_x)} \le
Ce^{-\de
t}\norm{f}_{W^{N_0+\abs{\al}+r-j}_1}\\
\intertext{and}
\normBig{\int_{U}e^{i(x\cdot\xi+\tau(\xi)t)}a_j(\xi)
\xi^{\al}\tau(\xi)^r \hat{f}(\xi)\,d\xi}_{L^2(\R^n_x)} \le
Ce^{-\de t}\norm{f}_{W^{\abs{\al}+r-j}_2}
\end{gather*}
for all $t>0$\textup{,} $N_0>n$\textup{,} multi-indices~$\al$,
$r\in\R$ and $f\in C_0^\infty(U)$.
\end{prop}

So, for all $t>0$,
\begin{gather*}
\normBig{D^r_tD_x^\al\Big(\int_{\R^n}e^{i(x\cdot\xi+\tau(\xi)t)}
a_j(\xi)\hat{f}(\xi)\,dx\Big)}_{L^\infty} \le Ce^{-\de
t}\norm{f}_{W_1^{N_1+\abs{\al}+r-j}}\,,\\
\normBig{D^r_tD_x^\al\Big(\int_{\R^n}e^{i(x\cdot\xi+\tau(\xi)t)}
a_j(\xi)\hat{f}(\xi)\,dx\Big)}_{L^2} \le Ce^{-\de
t}\norm{f}_{W_2^{\abs{\al}+r-j}}\,,
\end{gather*}
where $N_1>n$, $r\ge0$, $\al$ multi-index; interpolating then gives,
\begin{equation*}
\normBig{D^r_tD_x^\al\Big(\int_{\R^n}e^{i(x\cdot\xi+\tau(\xi)t)}
a_j(\xi)\hat{f}(\xi)\,dx\Big)}_{L^q} \le Ce^{-\de
t}\norm{f}_{W_p^{N_p+\abs{\al}+r-j}}\,,
\end{equation*}
where $p^{-1}+q^{-1}=1$, $1\leq p\le 2$, $N_p\ge
n\big(\frac{1}{p}-\frac{1}{q}\big)$, $r\ge0$,~$\al$ a multi-index
and $f\in C_0^\infty(\R^n)$. Thus, in this case we have
exponential decay of the solution.

\subsection{Bounded $\abs\xi$} 

Let us now consider the
terms of the first sum in \eqref{EQ:divideduprepforu}, the case of
low frequencies,
\begin{equation}\label{EQ:bddxiintegral}
\FT^{-1}(E_j\cutoffN\FT f)(x,t)=
\frac{1}{2\pi}\int_{\R^n}e^{ix\cdot\xi}\Big(\sum_{k=1}^m
e^{i\tau_k(\xi)t}
A_j^k(\xi,t)\Big)\cutoffN(\xi)\hat{f}(\xi)\,d\xi\,.
\end{equation}
Unlike in the case above, here the characteristic roots
$\tau_1(\xi),\dots,\tau_m(\xi)$ are not necessarily distinct at
all points in the support of the integrand (which is contained in
the ball of radius~$2N$ about the origin); in particular, this
means that the $A_j^k(\xi,t)$ genuinely depend on~$t$ and we have
no simple formula valid for them in the whole region.

For this reason, we begin by systematically separating
neighbourhoods of points where roots meet---referred to henceforth
as multiplicities---from the rest of the region, and then
considering the two cases separately.

First, we need to understand in what type of sets the roots
$\tau_k(\xi)$ can intersect:
\begin{lem}\label{LEM:multiplerootssetisnice}
The complement of the set of multiplicities of a linear strictly
hyperbolic constant coefficient partial differential
operator~$L(D_x,D_t)$\textup{,}
\begin{equation*}
S:=\set{\xi\in\R^n:\tau_j(\xi)\ne\tau_k(\xi)\text{ for all }j\ne
k}\,,
\end{equation*}
is dense in~$\R^n$.
\end{lem}
\begin{proof}
First note
\begin{equation*}
S=\set{\xi\in\R^n:\De_L(\xi)\ne0}\,,
\end{equation*}
where~$\De_L$ is the discriminant of~$L(\xi,\tau)$. Now,
by Sylvester's Formula (see, for example,~\cite{gelf+kapr+zele94}), $\De_L$
is a polynomial in the coefficients of
$L(\xi,\tau)$, which are themselves polynomials in~$\xi$.
Hence,~$\De_L$ is a polynomial in~$\xi$; as it is not identically
zero (for large~$\abs{\xi}$, the characteristic roots are distinct,
and hence it is non-zero at such points), it cannot be zero on an
open set, and hence its complement is dense in~$\R^n$.
\end{proof}

\begin{cor}\label{COR:multiplerootsetisnice}
Let $L(\xi,\tau)$ be a linear strictly hyperbolic constant
coefficient partial differential operator with characteristic
roots $\tau_1(\xi),\dots,\tau_m(\xi)$. Suppose
$\curlyM_{kl}\subset\R^n$ is a set such that
$\tau_{k}(\xi)=\tau_{l}(\xi)$, for some $k\ne l$, for all
$\xi\in\curlyM_{kl}$. For $\ep>0$, define
\begin{equation*}
\curlyM_{kl}^\ep:=\set{\xi\in\R^n:\dist(\xi,\curlyM_{kl})\leq\ep}\,;
\end{equation*}
denote the minimal $\nu\in\N$ such that
$\meas(\curlyM_{kl}^{\ep})\le C\ep^{\nu}$ for all sufficiently small
$\ep>0$ by $\codim\curlyM_{kl}$. Then $\codim\curlyM_{kl}\ge 1$.
\end{cor}
\begin{proof}
Follows straight from Lemma~\ref{LEM:multiplerootssetisnice}: the
fact that $\curlyM_{kl}$ has non-empty interior ensures that its
$\ep$-neighbourhood is bounded by~$C\ep$ in at least one dimension
for all small $\ep>0$.
\end{proof}

With this in mind, we subdivide integral~\eqref{EQ:bddxiintegral}: suppose~$L$
roots meet on a
set~$\curlyM$ with $\codim\curlyM=\ell$; without loss of
generality, assume the coinciding roots are
$\tau_1(\xi),\dots,\tau_L(\xi)$. By continuity, there exists an
$\ep>0$ such that only characteristic roots coinciding
with~$\tau_k(\xi)$, $k\in\set{1,\dots,L}$, in $\curlyM^{\ep}$
are~$\tau_1(\xi),\dots,\tau_L(\xi)$. Furthermore, we may assume
that $\pa\curlyM^{\ep}\in C^1$: for each $\ep>0$ there exists a
set~$S_\ep$ with~$C^1$ boundary such that $\curlyM^{\ep}\subset
S_\ep$ and $\meas(\curlyM^{\ep})\to \meas(S_\ep)$ as $\ep\to0$.
Then:
\begin{enumerate}[leftmargin=*]
\item Let $\cutoffM_{\curlyM,\ep}\in C^\infty(\R^n)$ be a smooth
function identically~$1$ on~$\curlyM^{\ep}$ and identically zero
outside $\curlyM^{2\ep}$; now consider the subdivision
of~\eqref{EQ:bddxiintegral}:
\begin{multline*}
\int_{B_{2N}(0)}e^{ix\cdot\xi}E_j(\xi,t)\hat{f}(\xi)\,d\xi
=\int_{B_{2N}(0)}e^{ix\cdot\xi}E_j(\xi,t)\cutoffM_{\curlyM,\ep}(\xi)
\hat{f}(\xi)\,d\xi\\+
\int_{B_{2N}(0)}e^{ix\cdot\xi}E_j(\xi,t)(1-\cutoffM_{\curlyM,\ep})(\xi)
\hat{f}(\xi)\,d\xi\,;
\end{multline*}
for the second integral, simply repeat the above procedure around
any root multiplicities in $B_{2N}(0)\setminus\curlyM^{\ep}$.
\item For the first integral, the case where the integrand is
supported on $\curlyM^{\ep}$, split off the coinciding roots from
the others:
\begin{multline}\label{EQ:intdivisionsplittingoffmults}
\int_{B_{2N}(0)}e^{ix\cdot\xi}E_j(\xi,t)\cutoffM_{\curlyM,\ep}(\xi)
\hat{f}(\xi)\,d\xi
\\=\int_{B_{2N}(0)}e^{ix\cdot\xi}\Big(\sum_{k=1}^L
e^{i\tau_k(\xi)t}A_j^k(\xi,t)\Big)
\cutoffM_{\curlyM,\ep}(\xi)\hat{f}(\xi)\,d\xi
\\+\int_{B_{2N}(0)}e^{ix\cdot\xi} \Big(\sum_{k=L+1}^m
e^{i\tau_k(\xi)t}A_j^k(\xi,t)\Big)
\cutoffM_{\curlyM,\ep}(\xi)\hat{f}(\xi)\,d\xi.
\end{multline}
\item For the first integral, we use techniques discussed in
Section~\ref{SEC:bddxiaroundmults} below to estimate it.

\item For the second there are two possibilities: firstly, two or
more of the roots $\tau_{L+1}(\xi),\dots,\tau_m(\xi)$ coincide in
$\curlyM^{2\ep}$---in this case, repeat the procedure above for
this integral. Alternatively, these roots are all distinct in
$\curlyM^{2\ep}$---in this case, it suffices to study each
integral separately as the $A_k^j(\xi,t)$ are independent of~$t$,
and thus the expression~\eqref{EQ:Ajkformula} is valid and we can
write
\begin{multline*}
\int_{B_{2N}(0)}e^{ix\cdot\xi} \Big(\sum_{k=L+1}^m
e^{i\tau_k(\xi)t}A_j^k(\xi,t)\Big)
\cutoffM_{\curlyM,\ep}(\xi)\hat{f}(\xi)\,d\xi\\=
\sum_{k=L+1}^m\int_{B_{2N}(0)}e^{i[x\cdot\xi+\tau_k(\xi)t]}
A_j^k(\xi)\cutoffM_{\curlyM,\ep}(\xi)\hat{f}(\xi)\,d\xi\;;
\end{multline*}
note that in this case we
may use that the region is bounded to ensure the continuous
functions are also bounded.
\end{enumerate}
Continue this procedure until all multiplicities are accounted for
in this way.

\subsubsection{Roots separated from the
real axis}\label{SEC:bddxiawayfromaxis}

The case where characteristic roots are separated from the real axis is similar
to that for large~$\abs{\xi}$. Let us assume that $\tau_k(\xi)$ has no multiplicities
in the set~$\Om$; now, a result similar to 
Proposition~\ref{PROP:rootsawayfromaxis} holds for general integrals of this
form, and thus
\begin{equation*}\label{EQ:estforrootawayfromaxisbddxi}
\normBig{D^r_tD_x^\al\Big(\int_{\Om}e^{i(x\cdot\xi+\tau(\xi)t)}
a(\xi)\hat{f}(\xi)\,dx\Big)}_{L^q} \le Ce^{-\de
t}\norm{f}_{L^p}\,,
\end{equation*}
where $p^{-1}+q^{-1}=1$, $1\leq p\le 2$, $N_p\ge
n\big(\frac{1}{p}-\frac{1}{q}\big)$, $r\ge0$,~$\al$ a multi-index,
$f\in C_0^\infty(\R^n)$, $\de>0$ is a constant such that
$\Im\tau(\xi)\ge\de$ for all $\xi\in\Om$ and $C\equiv
C_{\Om,r,\al,p}>0$. So, in this case we have also have exponential
decay of the solution.

\subsubsection{Roots meeting the real axis with finite
order}

In the case of bounded~$\abs{\xi}$, we must also consider the
situation where the phase function~$\tau(\xi)$ meets the real
axis. Suppose~$\xi_0\in\Om$ is such a point, i.e.\
$\Im\tau(\xi_0)=0$, while in each punctured ball around~$\xi_0$,
$B'_\ep(\xi_0)\subset\Om$, $\ep>0$, there exists $\xi\in
B'_\ep(\xi_0)$ so that $\Im\tau(\xi)>0$. Then, we claim
that~$\xi_0$ is a root of~$\Im\tau(\xi)$ of finite order~$s$:
indeed, if~$\xi_0$ were a zero of $\Im\tau(\xi)$ of infinite
order, then, by the analyticity of $\Im\tau(\xi)$ at~$\xi_0$
(which follows straight from the analyticity of~$\tau(\xi)$
at~$\xi_0$) it would be identically zero in a neighbourhood
of~$\xi_0$, contradicting the assumption.

In condition (H2) of Theorem \ref{THM:dissipative} we actually
have that there exist
constants $c_0,c_1>0$ such that, for all~$\xi$ sufficiently close
to~$\xi_0$,
\begin{equation*}
c_0\abs{\xi-\xi_0}^{s}\le\abs{\Im\tau(\xi)}\le
c_1\abs{\xi-\xi_0}^2\,.
\end{equation*}
Indeed, the Taylor expansion of $\Im\tau(\xi)$ around~$\xi_0$,
\begin{equation*}
\Im\tau(\xi)=\sum_{i=1}^n
\pa_{\xi_i}\Im\tau(\xi_0)(\xi_i-(\xi_0)_i)
+O(\abs{\xi-\xi_0}^2)\,,
\end{equation*}
is valid for~$\xi\in B_\ep(\xi_0)\subset\Om$ for some small
$\ep>0$. Now, if $\xi\in B_\ep(\xi_0)$, then $-\xi+2\xi_0\in
B_\ep(\xi_0)$ also. However,
\begin{equation*}
\Im\tau(-\xi+2\xi_0)=-\sum_{i=1}^n
\pa_{\xi_i}\Im\tau(\xi_0)(\xi_i-(\xi_0)_i)+O(\abs{\xi-\xi_0}^2)\,;
\end{equation*}
thus, for $\ep>0$ chosen small enough, this means that either
$\Im\tau(\xi)\le0$ or $\Im\tau(-\xi+2\xi_0)\le0$---contradicting
the hypothesis that $\Im\tau(\xi)\ge0$ for all $\xi\in\Om$; hence,
$\pa_{\xi_i}\Im\tau(\xi_0)=0$ for each $i=1,\dots,n$. In
conclusion, $\Im\tau(\xi)=O(\abs{\xi-\xi_0}^2)$ for all $\xi\in
B_\ep(\xi_0)$.

Now, we need the following result, which is based in the
calculation of the $L^p-L^q$ decay estimate for the dissipative
wave equation in~\cite{mats76}, but is here extended to a more
general situation so that it can be used on a wider class of
equations:
\begin{prop}\label{PROP:generalrootsmeetingaxis}
Let $\phi:U\to\R$\textup{,} $U\subset\R^n$ open\textup{,} be a
continuous function and suppose $\xi_0\in U$ such that
$\phi(\xi_0)=0$ and that $\phi(\xi)>0$ in a punctured open
neighbourhood of~$\xi_0$\textup{,} 
denoted by $V\setminus\set{\xi_0}$.
Furthermore\textup{,} assume that\textup{,} for some $s>0$,
there exists a constant $c_0>0$ such that\textup{,} for all
$\xi\in V$\textup{,}
\begin{equation*}
\phi(\xi)\ge c_0\abs{\xi-\xi_0}^{s}\,.
\end{equation*}
Then\textup{,} for any function~$a(\xi)$ that is bounded and
compactly supported in~$U$\textup{,} and for all $t\geq0$, $f\in
C_0^\infty(\R^n)$\textup{,} and $r\in\R$\textup{,}
\begin{gather}
\int_{V}e^{-\phi(\xi)t}\abs{\xi-\xi_0}^r\abs{a(\xi)}
\abs{\hat{f}(\xi)}\,d\xi \le
C\bract{t}^{-(n+r)/s}\norm{f}_{L^{1}}\,,\label{EQ:matsL1est}\\
\intertext{and} \normbig{e^{-\phi(\xi)t}\abs{\xi-\xi_0}^r a(\xi)
\hat{f}(\xi)}_{L^2(V)} \le
C\bract{t}^{-r/s}\norm{f}_{L^2}\,.\label{EQ:matsL2est}
\end{gather}
\end{prop}
\begin{proof}
First, we give a straightforward result that is useful in proving
each of the estimates:
\begin{lem}\label{LEM:simplelemma}
For each~$\rho,M\geq 0$ and $\varsigma,c>0$ there exists $C\equiv
C_{\rho,\varsigma,M,c}\ge0$ such that\textup{,} for all
$t\ge0$\textup{,}
\begin{gather*}
\int_0^{M}x^\rho e^{-cx^\varsigma t}\,dx\le
C\bract{t}^{-(\rho+1)/\varsigma} \text{ and } \sup_{0\le x\le
M}x^\rho e^{-cx^\varsigma t}\le C\bract{t}^{-\rho/\varsigma}\,.
\end{gather*}
\end{lem}
\begin{proof}
For $0\le t\le 1$, each is clearly bounded: the first by
$\frac{M^{\rho+1}}{\rho+1}$ and the second by~$M^\rho$. For $t>1$,
set $y=xt^{1/\varsigma}$; with this substitution, the first
becomes
\begin{equation*}
\int_0^{M t^{1/\varsigma}}y^\rho
t^{-\rho/\varsigma}e^{-cy^\varsigma}t^{-1/\varsigma}\,dy \le
t^{-(\rho+1)/\varsigma}\int_0^\infty y^\rho
e^{-cy^\varsigma}\,dy\,,
\end{equation*}
while the second becomes
\begin{equation*}
\sup_{0\le y\le Mt^{1/\varsigma}}y^\rho
t^{-\rho/\varsigma}e^{-cy^\varsigma}\le
t^{-\rho/\varsigma}\sup_{y\ge0}y^\rho e^{-cy^\varsigma}\,;
\end{equation*}
that the right-hand side of each is then bounded follows from
standard results.
\end{proof}
Returning to the proof of~\eqref{EQ:matsL1est}, as~$a(\xi)$ is
bounded in~$U$ by assumption, we have
\begin{equation*}
\int_{V}e^{-\phi(\xi)t}\abs{\xi-\xi_0}^r\abs{a(\xi)}
\abs{\hat{f}(\xi)}\,d\xi \le C\int_{V'}
e^{-\phi(\xi)t}\abs{\xi-\xi_0}^r \abs{\hat{f}(\xi)}\,d\xi\,,
\end{equation*}
where $V'=V\cap\supp a$; this, in turn, can be estimated in the
following manner using the hypothesis on $\phi(\xi)$ and
H\"older's inequality:
\begin{align*}
\int_{V'} e^{-\phi(\xi)t}\abs{\xi-\xi_0}^r
\abs{\hat{f}(\xi)}&\,d\xi \le C\int_{V'}
e^{-c_0\abs{\xi-\xi_0}^{s}t}\abs{\xi-\xi_0}^r
\abs{\hat{f}(\xi)}\,d\xi\\
&\le C\int_{V'}
e^{-c_0\abs{\xi-\xi_0}^{s}t}\abs{\xi-\xi_0}^r\,d\xi
\norm{\hat{f}}_{L^{\infty}(V')}\,.
\end{align*}
Then, transforming to polar coordinates and using the
Hausdorff--Young inequality, we find that, for some $\ep>0$
(chosen so that $V'\subset B_\ep(\xi_0)$, possible since~$a(\xi)$
is compactly supported),
\begin{multline*}
\int_{V'} e^{-c_0\abs{\xi-\xi_0}^{s}t}\abs{\xi-\xi_0}^r\,d\xi
\norm{\hat{f}}_{L^{\infty}(V')}\\ \le
C\int_{S^{n-1}}\int_0^\ep\abs{\eta}^{r+n-1}
e^{-c_0\abs{\eta}^{s}t}\,d\abs{\eta} d\omega\norm{f}_{L^{1}(\R^n)}\,,
\end{multline*}
Finally, by the first part of Lemma~\ref{LEM:simplelemma}, we find
\begin{align*}
\int_{V}e^{-\phi(\xi)t}\abs{\xi-\xi_0}^r\abs{a(\xi)}
\abs{\hat{f}(\xi)}\,d\xi &\le C\int_0^\ep y^{r+n-1}
e^{-c_0y^{s}t}\,dy\norm{f}_{L^{1}(\R^n)}\\&\le
C\bract{t}^{-(n+r)/s}\norm{f}_{L^{1}}\,.
\end{align*}
This completes the proof of the first part.

Now let us look at the second part. By the second part of
Lemma~\ref{LEM:simplelemma},
\begin{multline*}
\normbig{e^{-\phi(\xi)t}\abs{\xi-\xi_0}^ra(\xi)
\hat{f}(\xi)}_{L^2(V)}^2\le \int_{V'}
e^{-2c_0\abs{\xi-\xi_0}^{s}t}\abs{\xi-\xi_0}^{2r}
\abs{\hat{f}(\xi)}^2\,d\xi\\ \le
C\bract{t}^{-2r/s}\int_{V'}e^{-c_0\abs{\xi-\xi_0}^{s}t}
\abs{\hat{f}(\xi)}^2\,d\xi\,.
\end{multline*}
The H\"older inequality implies that
\begin{equation*}
\int_{V'} e^{-c_0\abs{\xi-\xi_0}^{s}t}
\abs{\hat{f}(\xi)}^2\,d\xi\le\sup_{V'}
\absbig{e^{-c_0\abs{\xi-\xi_0}^{s}t}}
\norm{\hat{f}}_{L^{2}(V')}^2\le C\norm{f}_{L^{2}}^2\,,
\end{equation*}
and together these give the required estimate~\eqref{EQ:matsL2est}.
\end{proof}

So, using this proposition, we have, for all $t>0$, and
sufficiently small $\ep>0$,
\begin{multline*}
\normBig{D^r_tD^\al_x\int_{B_\ep(\xi_0)}e^{i(x\cdot\xi+\tau(\xi)t)}
a(\xi) \hat{f}(\xi)\,d\xi}_{L^\infty(\R^n_x)}\\
\le\int_{B_\ep(\xi_0)}e^{-\Im\tau(\xi)t}\abs{a(\xi)}\abs{\tau(\xi)}^r
\abs{\xi}^{\al} \abs{\hat{f}(\xi)}\,d\xi \le
C\bract{t}^{-n/s}\norm{f}_{L^{1}}\,,
\end{multline*}
and, using the Plancherel Theorem,
\begin{multline*}
\normBig{D^r_tD^\al_x\int_{B_\ep(\xi_0)}e^{i(x\cdot\xi+\tau(\xi)t)}
a(\xi) \hat{f}(\xi)\,d\xi}_{L^2(\R^n_x)}\\ =
C\normbig{e^{i\tau(\xi)t}\tau(\xi)^r\xi^\al a(\xi)
\hat{f}(\xi)}_{L^2(B_\ep(\xi_0))} \le C\norm{f}_{L^2}\,;
\end{multline*}
here we have used that~$\abs{\xi}^{\abs{\al}}\abs{\tau(\xi)}^r\le
C$ for $\xi\in V^\prime$ for $r\in\N$, $\al$ a multi-index.

Thus, for all $t>0$,
\begin{equation*}\label{EQ:estforrootmeetingaxis}
\normBig{D^r_tD^\al_x\int_{B_\ep(\xi_0)}e^{i(x\cdot\xi+\tau(\xi)t)}a(\xi)
\hat{f}(\xi)\,d\xi}_{L^q(\R^n_x)} \le
C\bract{t}^{-\frac{n}{s}\big(\frac{1}{p}-\frac{1}{q}\big)}
\norm{f}_{L^{p}}\,,
\end{equation*}
where $1\leq p\le2$, $p^{-1}+q^{-1}=1$. 

\begin{rem}\label{rem:mats}
If $\xi_0=0$, then Proposition~\ref{PROP:generalrootsmeetingaxis}
further tells us that
\begin{equation*}
\normBig{D^r_tD^\al_x
\int_{B_\ep(0)}e^{i(x\cdot\xi+\tau(\xi)t)}a(\xi)
\hat{f}(\xi)\,d\xi}_{L^q(\R^n_x)}\le
C\bract{t}^{-\frac{n+\abs{\al}}{s}\big(\frac{1}{p}-\frac{1}{q}\big)}
\norm{f}_{L^{p}}\,.
\end{equation*}
If $\Re\tau(\xi_0)=0$, then 
under condition \eqref{EQ:imest} we have
$\abs{\tau(\xi)}\le
\abs{\Im\tau(\xi)}\le c_1\abs{\xi-\xi_0}^{s_1}$ for~$\xi$
near~$\xi_0$, and so we get
\begin{equation*}
\normBig{D^r_tD^\al_x
\int_{B_\ep(\xi_0)}e^{i(x\cdot\xi+\tau(\xi)t)}a(\xi)
\hat{f}(\xi)\,d\xi}_{L^q(\R^n_x)} \le
C\bract{t}^{-\big(\frac{n+ r s_1}{s}\big)\big(\frac{1}{p}-\frac{1}{q}\big)}
\norm{f}_{L^{p}}\,.
\end{equation*}
If both assumptions hold, we get the improvement from both 
cases, which is the estimate by
$C\bract{t}^{-\big(\frac{n+|\alpha|+r s_1}{s}\big)
\big(\frac{1}{p}-\frac{1}{q}\big)}.$
\end{rem}

\subsection{Estimates for bounded $\abs\xi$ around
multiplicities}\label{SEC:bddxiaroundmults}

Finally, let us turn to finding estimates for the first term
of~\eqref{EQ:intdivisionsplittingoffmults}, which we may write in
the form
\begin{equation*}
\int_{\Om}e^{ix\cdot\xi}\Big(\sum_{k=1}^L
e^{i\tau_k(\xi)t}A_j^k(\xi,t)\Big)
\cutoffM(\xi)\hat{f}(\xi)\,d\xi\,,
\end{equation*}
where the characteristic roots $\tau_1(\xi),\dots,\tau_L(\xi)$
coincide on a set~$\curlyM\subset\Om$ of codimension~$\ell$ (in
the sense of Corollary~\ref{COR:multiplerootsetisnice}),
$\Om\subset\R^n$ is a bounded open set and $\cutoffM\in
C_0^\infty(\Om)$.

Unlike in the case away from multiplicities of characteristic
roots, we have no explicit representation for the
coefficients~$A_j^k(\xi,t)$, which in turn means we cannot split
this into~$L$ separate integrals. To overcome this, we first show,
in Section~\ref{SEC:resofroots}, that a useful representation for
the above integral exists that allows us to use techniques
from earlier. Using this alternative representation, it is a
simple matter to find estimates in the case where the image of the
set~$\curlyM$ is separated from the real axis.
The argument may be extended to the case when it arises
on the real axis as a result of all the roots meeting the axis
with finite order. Such argument is more elaborate but
not necessary for Theorem \ref{THM:dissipative}.

\subsubsection{Resolution of multiple roots}\label{SEC:resofroots}
In this section, we find estimates for
\begin{equation*}
\sum_{k=1}^L e^{i\tau_k(\xi)t}A_j^k(\xi,t)\,,
\end{equation*}
where $\tau_1(\xi),\dots,\tau_L(\xi)$ coincide on a set~$\curlyM$
of codimension~$\ell$. For simplicity, first consider the simplest
case, $L=2$ and $\curlyM=\set{\xi_0}$; the general case works 
in a more involved but similar way. So, assume
\begin{equation*}
\tau_1(\xi_0)=\tau_2(\xi_0)\text{ and }
\tau_k(\xi_0)\ne\tau_1(\xi_0)\text{ for }k=3,\dots,m\,;
\end{equation*}
by continuity, there exists a ball of radius $\ep>0$
about~$\xi_0$, $B_\ep(\xi_0)$, in which the only root which
coincides with~$\tau_1(\xi)$ is~$\tau_2(\xi)$. Then:
\begin{lem}\label{LEM:2rootsmeetingresolution}
For all $t\ge0$ and $\xi\in B_\ep(\xi_0)$,
\begin{equation}\label{EQ:intersectingrootsbound}
\absBig{\sum_{k=1}^2e^{i\tau_k(\xi)t}A_j^k(\xi,t)}\le
C(1+t)e^{-\min(\Im\tau_1(\xi),\Im\tau_2(\xi))t}\,,
\end{equation}
where the minimum is taken over $\xi\in B_\ep(\xi_0)$.
\end{lem}
\begin{proof}
First, note that in the set
\begin{equation*}
S:=\{\xi\in\R^n:\tau_1(\xi)\ne\tau_k(\xi)\;\forall
k=2,\dots,m\text{ and}\,\tau_2(\xi)\ne\tau_l(\xi)\,\forall
l=3,\dots,m\}
\end{equation*}
the formula~\eqref{EQ:Ajkformula} is valid for $A_j^1(\xi)$ and
$A_j^2(\xi)$. Now, recall that \\
$E_j(\xi,t)=\sum_{k=1}^me^{i\tau_k(\xi)t}A^k_j(\xi,t)$ is the
solution to the Cauchy
problem~\eqref{EQ:ftcauchyprob},~\eqref{EQ:initialdataforEj}, and
thus is continuous; therefore, for all~$\eta\in\R^n$ such that
$\tau_1(\eta)\ne\tau_k(\eta)$ and $\tau_2(\eta)\ne\tau_k(\eta)$
for $k=3,\dots,m$ (but allow $\tau_1(\eta)=\tau_2(\eta)$),
we have
\begin{align*}
\sum_{k=1}^2e^{i\tau_k(\eta)t}A_j^k(t,\eta)
&=\lim_{\xi\to\eta}\big(e^{i\tau_1(\xi)t}A_j^1(\xi)+
e^{i\tau_2(\xi)t}A_j^2(\xi)\big)\,,
\end{align*}
provided~$\xi$ varies in the set~$S$ (thus, ensuring
$e^{i\tau_1(\xi)t}A_j^1(\xi)+ e^{i\tau_2(\xi)t}A_j^2(\xi)$ is
well-defined). Hence, to obtain~\eqref{EQ:intersectingrootsbound}
for all $\xi\in B_\ep(\xi_0)$, it suffices to show
\begin{equation*}
\absbig{e^{i\tau_1(\xi)t}A_j^1(\xi)+ e^{i\tau_2(\xi)t}A_j^2(\xi)}
\le Cte^{-\min(\Im\tau_1(\xi),\Im\tau_2(\xi))t}
\end{equation*}
for all $t\ge0$, $\xi\in
B'_\ep(\xi_0)=B_\ep(\xi_0)\setminus\set{\xi_0}$.

Now, for all $\xi\in B'_\ep(\xi_0)$, $t\ge0$,
\begin{multline}\label{EQ:splittingfor2roots}
e^{i\tau_1(\xi)t}A_j^1(\xi)+e^{i\tau_2(\xi)t}A_j^2(\xi)
\\=\sinh[(\tau_1(\xi)-\tau_2(\xi))t]
(e^{i\tau_2(\xi)t}A^1_j(\xi)-e^{i\tau_1(\xi)t}A^2_j(\xi))\\
+\cosh[(\tau_1(\xi)-\tau_2(\xi))t]
(e^{i\tau_2(\xi)t}A^1_j(\xi)+e^{i\tau_1(\xi)t}A^2_j(\xi))\,.
\end{multline}
Furthermore, we have the following estimates for all $\xi\in
B'_\ep(\xi_0)$, $t\ge0$:
\begin{multline}\label{EQ:sinhest}
\absbig{\sinh[(\tau_1(\xi)-\tau_2(\xi))t]
(A^1_j(\xi)e^{i\tau_2(\xi)t}- A^2_j(\xi)e^{i\tau_1(\xi)t})}\\ \le
Ct[\abs{e^{i\tau_2(\xi)t}}+\abs{e^{i\tau_1(\xi)t}}] \le
Cte^{-\min(\Im\tau_1(\xi),\Im\tau_2(\xi))t}\,,
\end{multline}
\begin{multline}\label{EQ:coshbound}
\abs{\cosh[(\tau_1(\xi)-\tau_2(\xi))t](A^1_j(\xi)e^{i\tau_2(\xi)t}
+A^2_j(\xi)e^{i\tau_1(\xi)t})}\\ \le
Cte^{-\min(\Im\tau_1(\xi),\Im\tau_2(\xi))t}\,.
\end{multline}
The proof of the first is simple: just note that
\begin{equation*}
\frac{\sinh[(\tau_1(\xi)-\tau_2(\xi))t]}{(\tau_1(\xi)-\tau_2(\xi))}
\to t\;\text{ as }\,(\tau_1(\xi)-\tau_2(\xi))\to0\,,
\end{equation*}
or, equivalently, as $\xi\to\xi_0$ through $S$, and
$A^k_j(\xi)(\tau_1(\xi)-\tau_2(\xi))$ is continuous in
$B_\ep(\xi_0)$ for $k=1,2$. The proof of the second is more technical
and uses the explicit representation~\eqref{EQ:Ajkformula} for the~$A_j^k(\xi)$
at points away from multiplicities of~$\tau_k(\xi)$;
otherwise it is similar and we omit it here.

Combining~\eqref{EQ:splittingfor2roots},~\eqref{EQ:sinhest} and
\eqref{EQ:coshbound} we have~\eqref{EQ:intersectingrootsbound},
which completes the proof of the lemma.

\end{proof}

Suppose now that
the characteristic roots $\tau_1(\xi),\dots,\tau_L(\xi)$,
$2\le L\le m$, coincide on a set $\curlyM$ of codimension~$\ell$,
and that $\tau_1(\xi)\ne\tau_k(\xi)$ for all $\xi\in\curlyM$ when
$k=L+1,\dots,m$. By continuity, we may take $\ep>0$ so that the
set $\curlyM^\ep=\set{\xi\in\R^n: \dist(\xi,\curlyM)\leq\ep}$
contains no points~$\eta$ at which
$\tau_1(\eta),\dots,\tau_L(\eta)=\tau_k(\eta)$ for
$k=L+1,\dots,m$. With this notation, we can extend
Lemma \ref{LEM:2rootsmeetingresolution} to the general
situation:
\begin{lem}\label{LEM:Lrootsmeetingest}
For all $t\ge0$ and $\xi\in\curlyM^\ep$\textup{,}
\begin{equation*}\label{EQ:multiplerootsrepbound}
\absBig{\sum_{k=1}^Le^{i\tau_k(\xi)t} A_j^k(\xi,t)} \le
C(1+t)^{L-1}e^{-t\min_{k=1,\dots,L}\Im\tau_k(\xi)}\,,
\end{equation*}
where the minimum is taken over $\xi\in B_\ep(\xi_0)$.
\end{lem}
Note that this estimate does not depend on the codimension
of~$\curlyM$ nor its geometric structure.

\subsubsection{Phase function separated from the real
axis}\label{SEC:phasefnsepfromaxis}

We now turn back to finding $L^p-L^q$ estimates for
\begin{equation*}
\int_{\Om}e^{ix\cdot\xi}\Big(\sum_{k=1}^L
e^{i\tau_k(\xi)t}A_j^k(\xi,t)\Big)
\cutoffM(\xi)\hat{f}(\xi)\,d\xi\,,
\end{equation*}
when $\tau_1(\xi),\dots,\tau_L(\xi)$ coincide on a set~$\curlyM$
of codimension~$\ell$; choose $\ep>0$ so that these roots do not
intersect with any of the roots
$\tau_{L+1}(\xi),\dots,\tau_m(\xi)$ in $\curlyM^\ep$.

Here we can assume that there exists $\de>0$ such that
$\Im\tau_k(\xi)\ge\de$ for all $\xi\in\curlyM^\ep$---so,
$\min_k\Im\tau_k(\xi)\ge\de$. For this, we use the same approach
as in Section~\ref{SEC:bddxiawayfromaxis}, but using
Lemma~\ref{LEM:Lrootsmeetingest} to estimate the sum. Firstly, the
$L^1-L^\infty$ estimate:
\begin{align*}
\normBig{D^r_tD_x^\al\Big(\int_{\Om}& e^{ix\cdot\xi}
\Big(\sum_{k=1}^L e^{i\tau_k(\xi)t}A_j^k(\xi,t)\Big)
\cutoffM(\xi)\hat{f}(\xi)\,dx\Big)}_{L^\infty(\R^n_x)}\\
&=\normBig{\int_{\Om}e^{ix\cdot\xi} \Big(\sum_{k=1}^L
e^{i\tau_k(\xi)t}A_j^k(\xi,t)\tau_k(\xi)^r\Big)
\xi^\al\cutoffM(\xi)\hat{f}(\xi)\,dx}_{L^\infty(\R^n_x)} \\
&\le \max_k\sup_\Om\abs{\tau_k(\xi)}^r
\int_{\Om}\absBig{\sum_{k=1}^L
e^{i\tau_k(\xi)t}A_j^k(\xi,t)}\abs{\xi}^{\abs{\al}}
\abs{\hat{f}(\xi)}\,dx\\
&\le C(1+t)^{L-1}e^{-\de t}\norm{\hat{f}}_{L^\infty(\Om)}\le
C(1+t)^{L-1}e^{-\de t}\norm{f}_{L^1}\,.
\end{align*}
Similarly, the $L^2-L^2$ estimate:
\begin{align*}
\normBig{D^r_tD_x^\al\Big(\int_{\Om}& e^{ix\cdot\xi}
\Big(\sum_{k=1}^L e^{i\tau_k(\xi)t}A_j^k(\xi,t)\Big)
\cutoffM(\xi)\hat{f}(\xi)\,dx\Big)}_{L^2(\R^n_x)}\\
&=\normBig{\Big(\sum_{k=1}^L
e^{i\tau_k(\xi)t}A_j^k(\xi,t)\tau_k(\xi)^r\Big)
\xi^\al\cutoffM(\xi)\hat{f}(\xi)}_{L^2(\Om)} \\
&\le C(1+t)^{L-1}e^{-\de t}\norm{\hat{f}}_{L^2(\Om)}\le
C(1+t)^{L-1}e^{-\de t}\norm{f}_{L^2}\,.
\end{align*}
Then,
\begin{multline*}
\normBig{D^r_tD_x^\al\Big(\int_{\Om} e^{ix\cdot\xi}
\Big(\sum_{k=1}^L e^{i\tau_k(\xi)t}A_j^k(\xi,t)\Big)
\cutoffM(\xi)\hat{f}(\xi)\,dx\Big)}_{L^q(\R^n_x)}\\ \le
C(1+t)^{L-1}e^{-\de t}\norm{f}_{L^p}\,,
\end{multline*}
where $p^{-1}+q^{-1}=1$, $1\le p\le2$. Once again, we have
exponential decay.

\section{Applications}\label{SEC:applications}

In this section we will briefly consider applications of Theorems
\ref{THM:overallmainthm} and \ref{THM:dissipative}
to Fokker-Planck equations and wave equations
with dissipation and negative mass. There are further
applications to Grad systems linearised near equilibrium
points where Theorems
\ref{THM:dissipative} and \ref{THM:dissspecial} immediately
yield the corresponding decay rates. The size of
these systems depends on the number of moments and dimension
of the space. Some examples of these systems and their
stability has been analysed in \cite{VR03}.

\subsection{Fokker--Planck Equation}

The classical Boltzmann equation for the particle 
distribution 
function $f=f(t,x,c)$, where $x,\mathbf{c}\in\R^n$, 
$n=1,2,3$, is 
\begin{equation*}
(\pa_t+\mathbf{c}\cdot\grad_x)f=S(f),
\end{equation*}
where $S(f)$ is the so-called integral of collisions.
The important special case of this equation is the
Fokker--Planck equation for the distribution of Brownian
particles, when the integral of collisions is linear and
is given by 
$$S(f)=\grad_\mathbf{c}\cdot(\mathbf{c}+\grad_{\mathbf{c}})f=
\sum_{k=1}^n \partial_{c_k}(c_k+\partial_{c_k})f.$$
In this case the kinetic Fokker-Planck equations takes the form
$$
\left(\pa_t+\sum_{k=1}^n c_k \partial_{x_k}\right)f(t,x,c)=
\sum_{k=1}^n \partial_{c_k}(c_k+\partial_{c_k})f.$$
The Hermite-Grad method of dealing with Fokker-Planck equation
consists in decomposing $f(t,x,\cdot)$ in the Hermite basis, i.e.
writing
$$f(t,x,c)=\sum_{|\alpha|\geq 0} \frac{1}{\alpha !}
m_\alpha(x,t) \psi^\alpha(c),$$ 
where
$\psi^\alpha(c)=(2\pi)^{-n/2}(-\partial_c)^\alpha
\exp(-\frac{|c|^2}{2})$ are Hermite functions. 
They are derivatives of the Maxwell distribution
$\psi^0$ which annihilates the integral of collisions
and form a complete orthonormal basis in the weighted
Hilbert space
$L^2_w(\R^n)$ with weight $w=1/\psi^0.$
This decomposition
yields the infinite system
$$
\pa_tm_\be(x,t)+\be_k\pa_{x_k}m_{\be-e_k}(x,t)+\\
\pa_{x_k}m_{\be+e_k}(x,t)+\abs{\be}m_\be(x,t)=0.
$$
The Galerkin approximation $f^N$ of the solution $f$ is
\begin{equation*}\mspace{-20mu}
f^N(t,x,c)=\sum_{0\le\abs{\al}\le N}\frac{1}{\al!}m_\al(x,t)
\psi^\alpha(c)\,,
\end{equation*}
with $m(x,t)=\{ m_\beta(x,t):\; 0\leq |\beta|\leq N\}$ being
the unknown function of coefficients.
For $m(x,t)$ one obtains the following system
of equations
\begin{equation*}
\pat m(x,t)+\sum_jA_j\paxj m(x,t)-iBm(x,t)=0,
\end{equation*}
where $B$ is a diagonal matrix, 
$B_{\alpha,\beta}=|\alpha|\delta_{\alpha,\beta},$
and the only non-zero elements of the matrix $A_j$ are
$a_j^{\alpha-e_j,\alpha}=\alpha_j$,
$a_j^{\alpha+e_j,\alpha}=1$. For details of these calculations
see \cite{vole+radk04}.
Hence, the dispersion equation for the system is
\begin{equation}
\begin{aligned}\label{EQ:FPpol}
&P(\tau,\xi)\equiv\det(\tau I+\sum_j A_j\xi_j-iB)=0\\
&P(\tau,0)=\det(\tau I-iB)=\tau\prod_{j=1}^{N}(\tau-ji)^
{\gamma_j}
\,. 
\end{aligned}
\end{equation}
Properties of this polynomial $P(\tau,\xi)$
have been extensively studied 
by Volevich and Radkevich in
\cite{vole+radk04}, who gave conditions and examples of
situations
when $\Im\tau_j(\xi)\ge 0$, for all $\xi\not=0$.
In our situation here we have to take additional care of
possible multiple roots, as is done in Theorem
\ref{THM:dissipative}.

Assume now that $P(\tau,\xi)$ is a stable polynomial, i.e.
its roots $\tau(\xi)$ satisfy $\Im\tau(\xi)\geq 0$ and
$\Im\tau(\xi)=0$ imply $\xi=0$. We will say that $P(\tau,\xi)$
is strongly stable if, moreover, its roots $\tau(\xi)$ satisfy
$\Im\tau(\xi)\geq\epsilon>0$ for large $\xi$.
It follows that we satisfy the conditions
of Theorem \ref{THM:dissipative} and we have to 
determine the 
order with which the characteristic arrives at the origin.
We have the following theorem about time decay of solutions
to Cauchy problems for equations with strongly
stable symbols.

\begin{thm}\label{THM:dissspecial}
Let a strongly stable polynomial $P(\tau,\xi)$ of 
order $m$ have a strictly 
hyperbolic principal part and assume that
$\partial_\tau P(0,0)\not=0$. Let $\alpha$ be
the multiindex of the smallest length $|\alpha|$ such that
$\partial_\xi^\alpha P(0,0)\not=0$.
Let $u(x,t)$ be the solution of the Cauchy problem
$P(D_t,D_x)u=0, \partial_t^l u|_{t=0}=f_l$,
$0\leq l\leq m-1$. 
Let $1\leq p\leq 2$, $2\leq q\leq\infty$, and $1/p+1/q=1$.
Then
$$ ||u(\cdot,t)||_{L^q}\leq C (1+t)^{-\frac{n}{|\alpha|}
\left(\frac1p-\frac1q\right)}
 \sum_{l=0}^{m-1} ||f_l||_{W_p^{N_l}},$$
where $N_l=\frac{2-p}{p}([n/2]+1)-l.$
Moreover, we have the estimate
\begin{equation*}
\normBig{\partial_t^r\partial^\beta_x
u(\cdot,t)}_{L^q(\R^n_x)}\le
C\bract{t}^{-\left(\frac{n+\abs{\beta}}{|\alpha|}+r\right)
\big(\frac{1}{p}-\frac{1}{q}\big)}
\sum_{l=0}^{m-1}
\norm{f_l}_{W^{N_l+\abs{\beta}+r}_p}
\,.
\end{equation*}
\end{thm}
Indeed, since the polynomial is strongly stable, the estimate
in Theorem \ref{THM:dissspecial} follows from Theorem
\ref{THM:dissipative} and Remark \ref{REM:dismore}.
The improvement in the last estimate for derivatives
comes from the fact that there is only one root $\tau$
such that $\tau(0)=0$ and so the last statement of
Remark \ref{REM:dismore} applies.
In certain cases it can be shown that actually 
$\abs{\al}=2$, in which case we have the same decay as for
dissipative wave equation.

Let us write the polynomial $P(\tau,\xi)$ from 
\eqref{EQ:FPpol} in the form 
$$P(\tau,\xi)=\sum_{j=0}^{M+1} (-i)^j P_j(\tau,\xi),$$
where $P_j$ is a homogeneous polynomial of order $M+1-j$,
and assume that $P(\tau,\xi)$ is strongly stable.
In \cite{vole+radk04} it was shown that
$P_{M+1}(\tau,\xi)\equiv 0$,
$\gamma_M P_M(\tau,\xi)=M!\tau$ for some $\gamma_M>0$,
and $\gamma_{M-1}P_{M-1}(\tau,\xi)=M!\sum_{k=2}^{M+1}
\frac{1}{k-1}\tau^2-M!|\xi|^2$ for some
$\gamma_{M-1}>0$. It can be readily verified now 
that conditions of Theorem \ref{THM:dissspecial} hold
with $|\alpha|=2$, from which we get the estimate with
$(1+t)^{-\frac{n}{2}\left(\frac1p-\frac1q\right)}.$

\subsection{Wave type equations with (negative) mass and dissipation}

\newcommand\mass{\mu}
\newcommand\diss{\delta}

Here we will show that we can still have time decay of
solutions if we allow the negative mass but exclude certain
low frequencies for Cauchy data. This is given in 
\eqref{EQ:wave} below. Nonnegative but
time dependent mass and
dissipation with oscillations have been considered before.
See, for example, \cite{HR03} and references therein.

Let us consider second order equations of the following form 
\begin{equation*}
\left\{
\begin{aligned} \pa_t^2u-c^2\lap u+\diss \pa_tu+\mass u=0\,,\\
u(0,x)=0,\;u_t(0,x)=g(x)\,.
\end{aligned}
\right.
\end{equation*}
Here $\diss$ is the dissipation and $\mass$ is the mass.
For simplicity, the first Cauchy data is taken to be zero.
The general case can be treated in the same way.
Let us now apply Theorem \ref{THM:overallmainthm} to the
analysis of this equation.
The associated characteristic polynomial is
\begin{equation*}
\tau^2-c^2\abs{\xi}^2-i\diss\tau-\mass=0\,,
\end{equation*}
which has roots
\begin{equation*}
\tau_\pm(\xi)=\frac{i\diss}{2}\pm\sqrt{c^2\abs{\xi}^2+\mass-\diss^2/4}\,.
\end{equation*}
Now, we have the following cases, which correspond to different
cases of Theorem \ref{THM:overallmainthm}:
\begin{itemize}
\item $\diss=\mass=0$. This is the wave equation.
\item $\diss=0$, $\mass>0$. This is the Klein--Gordon equation.
\item $\mass=0$, $\diss>0$. This is the dissipative wave equation.
\item $\diss<0$. In this case, 
$\Im\tau_-(\xi)\le\frac{\diss}{2}<0$ for
all $\xi$, hence we cannot expect any decay in general.
\item $\diss>0$, $\mass>0$. 
In this case the discriminant is always strictly
greater than $-\diss^2/4$, and thus the roots always lie in the 
upper half plane and are separated from the real axis.
So we have exponential decay.

Here is the main case for us:
\item $\diss\ge0$, $\mass<0$. In this case, note that
$\Im\tau_-(\xi)\ge0$ if and only if $c^2\abs{\xi}^2+\mass\ge0$,
i.e. the critical value is
$\abs{\xi}=\sqrt{\abs{\mass}}/{c}$. Therefore, the answer depends
on the Cauchy data $g$. In particular, if
$\supp \hat{g}$ is  contained in $c^2\abs{\xi}^2+\mass\geq 0$,
then we may get decay of some type. More precisely,
let $B(0,r)$ be the open ball with radius $r$ centred at the origin.
Then we have:
\begin{itemize}
\item if $g$ is such that $\supp\hat{g}\cap
B(0,\frac{\sqrt{\abs{\mass}}}{c})\ne\varnothing$, then we
have no decay;
\item if there is some
$\epsilon>0$ such that $\supp\hat{g}\subset \R^n\setminus
B(0,\frac{\sqrt{\abs{\mass}}}{c}+\ep)$, then the
roots are either separated from the real axis (if $\de>0$), and we
get exponential decay, or lie on the real axis 
(if $\de=0$), and we
get Klein--Gordon type behaviour 
(since the Hessian of $\tau$ is nonsingular).
\item if, for all $g$, $\supp\hat{g}\subset\R^n\setminus
B(0,\frac{\sqrt{\abs{\mass}}}{c})
=\left\{|\xi|\geq \frac{\sqrt{\abs{\mass}}}{c}\right\}$ , 
then again we must consider
$\de=0$ and $\de>0$ separately.

If $\de=0$, then the roots lie completely on the real axis, 
and they meet on the sphere
$|\xi|=\sqrt{|\mu|}/{c}$. It follows from \eqref{EQ:around}
with $L=2$ and $\ell=1$ that, although the representation
of solution as a sum of Fourier integrals breaks down at the
sphere, the solution is still bounded
in a ($1/t$)-neighbourhood of the sphere. In its complement we
can get the decay.

If $\de>0$, then the root
$\tau_{-}$ comes to
the real axis at $\abs{\xi}=\frac{\sqrt{\abs{\mass}}}{c}$, in
which case we get the decay
\begin{equation}\label{EQ:wave}
||u(\cdot,t)||_{L^q}\leq C (1+t)^{-\left(\frac1p-\frac1q
\right)} ||g||_{L^p}.
\end{equation}
Indeed, in this case
the order of the root $\tau_{-}$ at the axis is one,
i.e. estimate \eqref{EQ:disttau} holds with $s=1$.
Here $1/p+1/q=1$ and $1\leq p\leq 2$. 
Note also that compared to the case of 
no mass when $\ell=n$, the codimension of the set
$\set{\xi\in\R^n:\abs{\xi}=\frac{\sqrt{\abs{\mass}}}{c}}$ is
$\ell=1$. We can apply the last case of Part II of
Theorem \ref{THM:overallmainthm} with 
$L=1$ and $s=\ell=1$ which 
gives estimate \eqref{EQ:wave}.
\end{itemize}

\end{itemize}


\end{document}